\documentclass[a4paper,12pt]{amsart}
\usepackage{amsopn,amsfonts,amssymb,amsthm,mathtools,amsmath}
\usepackage{color}
\usepackage{latexsym}

\usepackage[T1]{fontenc}
\usepackage{lmodern}
\usepackage{textcomp}
\usepackage{hyperref}
\hypersetup{colorlinks=true,linkcolor=blue,citecolor=blue}
\usepackage[utf8]{inputenc}
\usepackage{enumerate}
\usepackage[shortlabels]{enumitem}
\setlist[enumerate, 1]{\sc(1)}
\usepackage{float}
\usepackage{anysize}
\marginsize{2cm}{2cm}{2cm}{2cm}
\usepackage[justification=centering]{caption}

\usepackage{enumerate}

\def\u{\mathfrak{u}}
\def\b{\mathfrak{b}}

\def\g{\mathfrak{g}}

\def\e{\mathfrak{e}}
\def\a{\mathfrak{a}}

\def\sp{\mathfrak{sp}}

\def\u{\mathfrak{u}}
\def\gl{\mathfrak{gl}}

\def\X{\mathfrak{X}}

\def\B{\mathcal{B}}

\def\W{\mathcal{W}}
\def\l{\ell}

\def\C{\mathbb{C}}
\def\R{\mathbb{R}}

\def\Z{\mathbb{Z}}
\def\N{\mathbb{N}}

\newcommand{\Sym}{{\operatorname{Sym}}}

\def\ad{\operatorname{ad}}

\def\Tr{\operatorname{Tr}}

\def\I{\operatorname{I}}

\newcommand{\Ker}{\operatorname{Ker}}

\DeclareMathOperator{\im}{Im}

\def\alt{\raise1pt\hbox{$\bigwedge$}}
\def\pint{\langle \cdotp,\cdotp \rangle }
\def\la{\langle}
\def\ra{\rangle}

\theoremstyle{plain}
\newtheorem{theorem}{\bf Theorem}[section]
\newtheorem{corollary}[theorem]{\bf Corollary}
\newtheorem{proposition}[theorem]{\bf Proposition}
\newtheorem{lemma}[theorem]{\bf Lemma}

\theoremstyle{definition}

\newtheorem{example}[theorem]{\bf Example}

\theoremstyle{remark}
\newtheorem{remark}[theorem]{\bf Remark}
\newcommand{\ri}{{\rm (i)}}
\newcommand{\rii}{{\rm (ii)}}

\newcommand{\matriz}[1]{\ensuremath{\begin{bmatrix}#1\end{bmatrix}}}
\newcommand{\comillas}[1]{``#1''}

\title{Harmonic almost complex structures}

\author{Adri\'an Andrada}
\email{adrian.andrada@unc.edu.ar}

\author{Alejandro Tolcachier}
\email{atolcachier@unc.edu.ar}

\address{FAMAF, Universidad Nacional de C\'ordoba and CIEM-CONICET, Av. Medina Allende s/n, X5000HUA C\'ordoba, Argentina}

\thanks{This work was partially supported by CONICET, SECyT-UNC and FONCyT (Argentina) and the MATHAMSUD Regional Program 21-MATH-06.}

\subjclass[2020]{53C15, 53C55, 22E25, 22E40}
\keywords{Harmonic almost complex structure, Gray-Hervella classes, almost abelian Lie group, lattice, solvmanifold}

\title[Harmonic almost complex structures]{Harmonic almost complex structures on almost abelian Lie groups and solvmanifolds}

\begin{document}
	
	\begin{abstract}
		An almost abelian Lie group is a solvable Lie group with a codimension one normal abelian subgroup. We characterize almost Hermitian structures on  almost abelian Lie groups where the almost complex structure is harmonic with respect to the Hermitian metric. Also, we adapt the Gray-Hervella classification of almost Hermitian structures to the family of almost abelian Lie groups. We provide several examples of harmonic almost complex structures in different Gray-Hervella classes on some associated compact almost abelian solvmanifolds. 
	\end{abstract}
	
	\maketitle
	
	\section{Introduction}
	
	Let $(M,g)$ be a $2n$-dimensional connected Riemannian manifold. An almost complex structure on $M$ is a $(1,1)$-tensor $J$ that satisfies $J_p^2=-\I_{T_pM}$ for all $p\in M$, and, moreover, if $J$ is orthogonal (i.e.,  $g(J\cdot,J\cdot)=g(\cdot,\cdot)$) then the pair $(J,g)$ is called an almost Hermitian structure and $(M,J,g)$ is called an almost Hermitian manifold. Note that, due to the condition $J^2=-\operatorname{I}$, $J$ is orthogonal if and only if $J$ is skew-symmetric.
	
	When $J$ is integrable, i.e. $J$ is the almost complex structure underlying a holomorphic atlas on $M$ making it a complex manifold, then $J$ is called a complex structure and the pair $(J,g)$ is known simply as a Hermitian structure. We recall that the integrability of the almost complex structure $J$ is equivalent to the vanishing of the Nijenhuis tensor associated to $J$, that is, the $(1,2)$-tensor $N$ defined by 
	\begin{equation}\label{eq:N}
		N(X,Y)=[X,Y]+J([JX,Y]+[X,JY])-[JX,JY], \quad X,Y\in \mathfrak{X}(M).
	\end{equation}
	
	If $(M,g)$ admits an orthogonal almost complex structure then it admits many of them, so that it would be important to be able to distinguish the ``best'', or the ``optimal'', compatible almost complex structure. Due to its rich geometric properties, it is clear that the best almost complex structures are the K\"ahler ones, that is, almost complex structures $J$ satisfying $\nabla J=0$, where $\nabla$ denotes the Levi-Civita connection of $g$. This condition is equivalent to the integrability of $J$ together with $d\omega=0$, where $\omega$ is the fundamental $2$-form defined by $\omega=g(J\cdot,\cdot)$. Nevertheless, due to strong topological restrictions, there are many compact almost Hermitian manifolds that do not admit K\"ahler structures. 
	
	One possible approach to detect an optimal almost complex structure is to consider functionals arising from the given geometrical data and then look for extrema of those functionals. A natural functional to analyze is the Dirichlet energy functional $E$ defined on the set of orthogonal almost complex structures on $M$, that is,
	\begin{equation}\label{eq:energy}
		E(J):=\int_M \|\nabla J\|^2 \operatorname{vol}_g, 
	\end{equation}
	when $M$ is compact. 
	
	The critical points of this functional are called \textit{harmonic} almost complex structures. It was shown in \cite{Wood-Crelle} (see also \cite{Wood}) that an orthogonal almost complex structure $J$ is harmonic if and only if the corresponding Euler-Lagrange equation is satisfied:
	\begin{equation}\label{eq:harmonic J}
		[J,\nabla^*\nabla J]=0, 
	\end{equation}
	where $\nabla^*\nabla J$ is the \textit{rough Laplacian} of $J$ defined by $\nabla^*\nabla J=\Tr \nabla^2J$, that is, the operator 
	\begin{equation}\label{eq:rough}
		(\nabla^*\nabla J)(X)= \sum_{i=1}^{2n} (\nabla^2_{e_i,e_i} J)(X), \quad X\in\mathfrak{X}(M),
	\end{equation} 
	where $\{e_1,\ldots,e_{2n}\}$ is a local orthonormal frame on $M$, and the second covariant derivative of $J$ is $(\nabla^2_{U,V} J)(W)=(\nabla_{U} (\nabla_{V} J))(W)-(\nabla_{\nabla_{U} V} J)(W)$. It is clear that $\nabla^*\nabla J$ is a $(1,1)$-tensor on $M$. On a non-compact manifold, one can either take the Euler-Lagrange equation $[J,\nabla^*\nabla J]=0$ as the definition of harmonicity or define the energy on an open subset with compact closure of the manifold and consider variations with compact support included in this subset, and the resulting equation on the open subset is the same as in the compact case.
	
	It follows clearly from \eqref{eq:energy} that K\"ahler structures are harmonic; indeed, they are absolute minimizers. In the more general almost Hermitian setting, we recall that almost Hermitian structures were divided into classes by Gray and Hervella in \cite{GH}, according to properties of the associated tensor $\nabla\omega$ (see \S\ref{section:GH} below). The almost complex structures in some of these classes were shown to be always harmonic; for instance, nearly K\"ahler manifolds in \cite{Wood} (such as $S^6$ with the almost complex structure induced by octonionic multiplication), and balanced or locally conformally K\"ahler manifolds in \cite{GM}. Other examples of harmonic complex structures appear on Calabi-Eckmann manifolds, which, as Riemannian manifolds, are the product of two odd-dimensional spheres equipped with the round metrics \cite{Wood}. Harmonic almost complex structures on $4$-dimensional Riemannian manifolds were studied in \cite{DM1,DM2,DUM}. More recently, He and Li introduced in \cite{HL} the \textit{harmonic heat flow} for almost complex structures compatible with a fixed Riemannian metric, which is  a tensor-valued version of the  harmonic map heat equation first studied by Eells-Sampson \cite{ES}. Generalizing the case of Calabi-Eckmann manifolds we showed in \cite{AT} that in a 2-parameter family of Hermitian structures $(J_{a,b},g_{a,b})$ on a product of Sasakian manifolds, the integrable almost complex structure $J_{a,b}$ is harmonic with respect to $g_{a,b}$.
	
	The Euler-Lagrange equation $[J,\nabla^*\nabla J]=0$ is difficult to verify on a general almost Hermitian manifold. Therefore we will deal with this problem in the case of a special family of compact manifolds, namely almost Hermitian solvmanifolds. We recall that a solvmanifold is a compact quotient $M_\Gamma=\Gamma\backslash G$  of a simply connected solvable Lie group $G$ by a discrete subgroup $\Gamma$.  They constitute a fruitful and interesting source of examples in (almost) Hermitian, symplectic and $G_2$ geometry, among others (see for instance \cite{FG,FKV,KL,MSa,UV}). Indeed, when studying the properties of geometric structures on $M_\Gamma$ induced by a left invariant one on the Lie group $G$, we can work at the Lie algebra level and consider the associated linear objects. In our case, on such a solvmanifold we will consider almost Hermitian structures induced by a left invariant almost Hermitian structure on the Lie group $G$. Consequently, it will be enough to verify equation \eqref{eq:harmonic J} at the Lie algebra level and thus, the harmonicity condition becomes an algebraic problem on a Lie algebra (which is not necessarily easy to check).
	
	Our aim is to study the harmonicity of invariant almost complex structures in a particular family of almost Hermitian solvmanifolds. Indeed, we will restrict ourselves to the case of \textit{almost abelian solvmanifolds} equipped with an invariant almost Hermitian structure. A solvmanifold $\Gamma\backslash G$ is called almost abelian when $G$ is an almost abelian Lie group, that is, its Lie algebra $\g$ has a codimension one abelian ideal. This class of solvmanifolds has been intensely studied lately, and results about their (almost) Hermitian geometry can be found in \cite{AO,AL,FP-1,FP,LRV}, to name a few. 
	
	Any $2n$-dimensional almost abelian Lie algebra is completely determined by a $(2n-1)\times (2n-1)$ real matrix $A$ and we interpret the harmonicity condition in terms of $A$. Using this we provide several examples of harmonic almost complex structures in different Gray-Hervella classes both in dimension 4 and in dimensions $2n\geq 6$. In particular, we recover the well-known almost Kähler structure on the Kodaira-Thurston manifold given by Abbena in \cite{Abbena}. This is an example of an almost abelian nilmanifold equipped with an invariant almost Hermitian structure, whose almost complex structure is harmonic according to \cite{Wood}.
	
	The outline of the article is the following. We begin with some preliminaries about solvmanifolds in \S\ref{section:prelim}. In \S\ref{section:main} we characterize the almost abelian Lie algebras equipped with an almost Hermitian structure with harmonic almost complex structure (Theorem \ref{theorem:main}). In \S\ref{section:GH} we describe in a unified way the Gray-Hervella classes on almost abelian Lie algebras in order to study the interplay between the harmonicity condition and almost Hermitian geometry. We carry out this study in \S \ref{section:examples}, which is completely devoted to giving examples of different phenomena combining harmonicity and the Gray-Hervella classes. Finally, we also analyze the relation between harmonicity and SKT geometry, which is a special class of Hermitian structures defined by $\partial\bar{\partial}\omega=0$ that does not appear in Gray-Hervella classification.
	
	\medskip 
	
	\textbf{Acknowledgments.} 
	The authors are grateful to Eric Loubeau, Andrés Moreno, Andrei Moroianu and Henrique Sá Earp for their useful comments and suggestions as well as to the reviewer for his/her careful reading of the manuscript. The authors would also like to thank the hospitality of the Instituto de Matemática, Estatística e Computação Científica at UNICAMP (Brazil), where this work was initiated. 
	
	\medskip
	
	\section{Preliminaries}\label{section:prelim}
	
	In this section we recall some facts about solvmanifolds and their invariant geometric structures. 
	
	A discrete subgroup $\Gamma$ of a simply connected Lie group $G$ is called a \textit{lattice} if the quotient $\Gamma\backslash G$ is compact. According to \cite{Mi}, if such a lattice exists then the Lie group must be unimodular\footnote{A connected Lie group is called \textit{unimodular} if $\Tr \ad_x=0$ for all $x$ in its Lie algebra.}. The quotient $\Gamma\backslash G$ is called a solvmanifold if $G$ is solvable and a nilmanifold if $G$ is nilpotent, and it is well known that $\pi_1(\Gamma\backslash G)\cong \Gamma$. Moreover, the diffeomorphism class of solvmanifolds is determined by the isomorphism class of the corresponding lattices, as the following result shows:
	
	\begin{theorem}\label{theorem:mostow}\cite{Mo} 
		Two solvmanifolds with isomorphic fundamental groups are diffeomorphic.
	\end{theorem}
	
	\smallskip
	
	It is not easy to determine whether a given unimodular solvable Lie group admits a lattice. However, there is such a criterion for nilpotent Lie groups. Indeed, it was proved by Malcev in \cite{Ma} that a nilpotent Lie group admits a lattice if and only if its Lie algebra has a rational form, i.e. there exists a basis of the Lie algebra such that the corresponding structure constants are all rational. More recently, in \cite{B}, the existence of lattices in simply connected solvable Lie groups up to dimension $6$ was studied. A general result proved by Witte in \cite[Proposition 8.7]{Wi} states that only countably many non-isomorphic simply connected Lie groups admit lattices. 
	
	\ 
	
	Let $G$ be a simply connected solvable Lie group with a left invariant almost Hermitian structure $(J,g)$, i.e. the left translations $L_x:G\to G$ defined by $L_x(h)=xh$ for $h\in G$ are both almost complex maps and isometries. 
	The left invariant almost Hermitian structure $(J,g)$ is uniquely determined by its values at the identity, therefore we can work at the Lie algebra level. So, if $\g$ denotes the Lie algebra of $G$, then $J$ induces an endomorphism of $\g$, still denoted $J:\g\to\g$, and $g$ induces an inner product $\pint$ on $\g$, such that 
	\[ J^2=-\operatorname{I}, \qquad \langle J\cdot,J\cdot\rangle=\pint. \] 
	These two conditions imply: $\langle J\cdot,\cdot \rangle=-\langle \cdot, J\rangle$. The pair $(J,\pint)$ will be called an almost Hermitian structure on $\g$. Moreover, the almost complex structure $J$ on $G$ is integrable if and only if $J:\g\to\g$ satisfies \eqref{eq:N} when $X,Y$ are elements in $\g$.
	
	The Levi-Civita connection $\nabla$ of the left invariant metric $g$ satisfies $\nabla_XY\in\g$ whenever $X,Y\in\g$. Therefore, $\nabla$ is uniquely determined by the bilinear map $\nabla:\g\times \g\to\g$, $(X,Y)\mapsto \nabla_XY$, defined by the Koszul formula
	\[ 2\langle \nabla_XY,Z\rangle = \langle [X,Y],Z\rangle -\langle [Y,Z],X\rangle +\langle [Z,X],Y\rangle .\]
	
	As a consequence, the rough Laplacian $\nabla^*\nabla J$ (as in \eqref{eq:rough}) associated to the left invariant almost Hermitian structure $(J,g)$ on $G$ is also a left invariant tensor, thus it is determined by its restriction to $\g$, and therefore it defines an endomorphism $\nabla^*\nabla J:\g\to\g$. It is clear that $[J,\nabla^*\nabla J]=0$ on $G$ if and only if $[J,\nabla^*\nabla J]=0$ on $\g$, since both $J$ and $\nabla^*\nabla J$ are left invariant tensors. Hence, in order to determine if the almost complex structure $J$ is harmonic with respect to $g$, it is enough to verify that $[J,\nabla^*\nabla J]$ vanishes on left invariant vector fields. Hence, we will say that an orthogonal almost complex structure $J$ on $(\g,\pint)$ is harmonic if 
	\[ [J,\nabla^*\nabla J]=0\quad \text{on } \g.\] 
	
	If $(J,g)$ is a left invariant almost Hermitian structure on the simply connected solvable Lie group $G$ and $G$ admits a lattice $\Gamma$ then the solvmanifold $\Gamma\backslash G$ inherits a unique almost Hermitian structure, still denoted by $(J,g)$ and called invariant, such that the projection $\pi:G\to \Gamma\backslash G$ is locally a holomorphic isometry. It is easy to verify that the invariant almost complex structure $J$ is harmonic on $(\Gamma\backslash G,g)$ if and only if $J$ is harmonic on $(\g,\pint)$. 
	
	To sum up, when considering the case of a Lie group equipped with a left invariant almost Hermitian structure, or a solvmanifold equipped with an invariant almost Hermitian structure, the harmonicity equation \eqref{eq:harmonic J} reduces to an algebraic problem on a Lie algebra, and this makes it easier to solve, even though it could still be involved. 
	
	\medskip 
	
	\subsection{Almost abelian Lie groups}
	
	A $d$-dimensional Lie group $G$ is said to be {\em almost abelian} if its Lie algebra $\g$ has a codimension one abelian ideal. Such a Lie algebra will be called almost abelian, and it can be written as $\g= \R e_0 \ltimes \mathfrak{u}$, where $\mathfrak u$ is an abelian ideal of $\g$, and $\R$ is generated by some $e_0\in \g$. After choosing a basis of $\u$, we may identify $\u$ with the abelian Lie algebra $\R^{d-1}$ and we may write $\g=\R e_0\ltimes_L \R^{d-1}$ for some $L\in \mathfrak{gl}(d-1,\R)$. 
	
	Accordingly, the Lie group $G$ can be written as a semidirect product $G=\R\ltimes_\phi \R^{d-1}$, where the action is given by $\phi(t)=e^{tL}$. We point out that a non-abelian almost abelian Lie group is $2$-step solvable, and it is nilpotent if and only if the operator $L$ is nilpotent. 
	
	Concerning the isomorphism classes of almost abelian Lie algebras, there is the following result (see for instance \cite[Proposition 1]{Fr}):
	
	\begin{lemma}\label{lemma:isomorphic}
		Two almost abelian Lie algebras $\g_1=\R \ltimes_{L_1} \R^{d-1}$ and 
		$\g_2=\R \ltimes_{L_2}\R^{d-1}$ are isomorphic if and only if there exists $c\neq 0$ such that
		$cL_1$ and $L_2$ are conjugate in $\mathfrak{gl}(d-1,\R)$. 
	\end{lemma}
	
	\medskip
	
	An important feature of almost abelian Lie groups is that there exists a criterion to determine whether they admit lattices or not.
	
	\begin{proposition}\label{proposition:bock}\cite{B}
		Let $G=\R\ltimes_\phi\R^m$ be an almost abelian Lie group. Then $G$ admits a lattice $\Gamma$ if and only if there exists $t_0\neq 0$ such that $\phi(t_0)$ is conjugate to an integer matrix in $\operatorname{SL}(m,\Z)$. In this case, the lattice $\Gamma$ is given by $\Gamma=t_0\Z \ltimes_\phi P\Z^m$, where $E=P^{-1}{\phi(t_0)}P\in \operatorname{SL}(m,\Z)$.
	\end{proposition}
	
	\smallskip
	
	\begin{remark}\label{remark:gamma-isomorphic}
		Note that $\Gamma$ is isomorphic as an abstract group to $\Z\ltimes_E \Z^m$, where the multiplication is given by \[(k,(r_1,\ldots,r_m))\cdot (\l,(s_1,\ldots,s_m))=\left(k+\l,(r_1,\ldots,r_m)+E^k (s_1,\ldots,s_m)\right).\]
	\end{remark}
	
	\medskip
	
	\section{Harmonic almost complex structures on almost abelian Lie algebras}\label{section:main}
	
	Let $G$ be a $2n$-dimensional almost abelian Lie group equipped with a left invariant almost Hermitian structure $(J,g)$, determined by an endomorphism $J:\g\to\g$ and an inner product $\pint$ on $\g$, where $\g$ denotes the Lie algebra of $G$.
	
	If $\u$ denotes the codimension one abelian ideal in $\g$, let us choose $e_0\in \g$ with $|e_0|=1$ orthogonal to $\u$, so that we can decompose $\g$ orthogonally as 
	\[\g=\R e_0\ltimes \u=\R e_0\ltimes_L \R^{2n-1}\] 
	for some $L\in \mathfrak{gl}(2n-1,\R)$. Let us denote by $\a$ the maximal $J$-invariant subspace of $\u$, that is, $\a=\u\cap J\u$. It is clear that $\dim \a= 2n-2$. Let $\{e_0, e_1,\ldots,e_{2n-1}\}$ be an orthonormal basis of $\g$ compatible with $J$, i.e. $Je_{2i}=e_{2i+1}$, $0\leq i\leq n-1$. Note that $\u=\text{span}\{e_1,\ldots,e_{2n-1}\}$ and $\a=\text{span}\{e_2,\ldots,e_{2n-1}\}$. We will denote $J':=J|_\a:\a\to\a$.
	
	We decompose the endomorphism $L$ of $\u$ according to the orthogonal decomposition $\u=\R e_1\oplus \a$:
	\begin{equation}\label{eq:matrix L} 
		L=\left[
		\begin{array}{c|ccc}       
			\mu & & w_0^t &\\
			\hline
			& & & \\
			v_0 &  & D & \\
			& & & 
		\end{array}
		\right] 
	\end{equation} 
	for some $\mu\in\R$, $v_0,w_0\in \a$ and $D\in\mathfrak{gl}(2n-2,\R)$. 
	
	\smallskip
	
	\begin{remark}\label{remark:change}
		If $w_0=0$ and $0\neq v_0\in \im(D-\mu \I)$ then we can perform a holomorphic change of basis such that $v_0=0$. Indeed, if $v_0=Dx-\mu x$ for some $x\in \a$, then  let us set $e_0'=e_0+Jx$, $e_1'=Je_0=e_1-x$. Now 
		\[ [e_0',e_1']=\mu e_1+v_0-Dx=\mu e_1', \qquad \ad_{e_0'}|_\a=D.\] 
		Note that this change of basis does not preserve the inner product, since the new basis $\{e_0',e_1', e_2,\ldots, e_{2n-1}\}$ is not orthonormal.
	\end{remark}
	
	\medskip 
	
	We will need the decomposition of $L$ into its symmetric and skew-symmetric components $L=S+A$, with 
	\begin{equation}\label{eq: desc. L}
		S=\left[
		\begin{array}{c|ccc}       
			\mu & & \gamma^t &\\
			\hline
			& & & \\
			\gamma &  & D_s & \\
			& & & 
		\end{array}\right] 
		, \quad A= \left[
		\begin{array}{c|ccc}       
			0 & & -\rho^t &\\
			\hline
			& & & \\
			\rho &  & D_a & \\
			& & & 
		\end{array}
		\right], 
	\end{equation}
	
	Here, \[\gamma=\frac12(v_0+w_0), \quad \rho=\frac12(v_0-w_0),\] and $D_s,\,D_a$ are the symmetric and skew-symmetric components of $D$, respectively. Since $L=\ad_{e_0}|_\u$, we may also set $Le_0=0$, and this implies $Se_0=Ae_0=0$.
	
	\medskip
	
	In order to determine whether the almost complex structure $J$ on $(\g,\pint)$ is harmonic we will use the following lemma, which holds for any almost Hermitian manifold.
	
	\begin{lemma}
		Let $(M,J,g)$ be a $2n$-dimensional almost Hermitian manifold and $\{e_i\}_{i=1}^{2n}$ be an orthonormal frame on an open subset of $M$. Then 
		\begin{equation}\label{eq:conm-explicit}
			[J,\nabla^* \nabla J](X)=2 \sum_{i=1}^{2n} \left( \nabla_{e_i} J \nabla_{e_i} J X-J\nabla_{e_i} J\nabla_{e_i} X -J(\nabla_{\nabla_{e_i} e_i} J)X\right)
		\end{equation}
	\end{lemma}
	\begin{proof}
		Let $X\in\X(M)$, then
		\begin{align*} 
			(\nabla^*\nabla J)(X)&=\sum_{i=1}^{2n}  (\nabla_{e_i} (\nabla_{e_i} J))(X)-(\nabla_{\nabla_{e_i} e_i} J)(X)\\
			&=\sum_{i=1}^{2n}  \nabla_{e_i}  ((\nabla_{e_i} J)(X))-(\nabla_{e_i} J)(\nabla_{e_i} X)-(\nabla_{\nabla_{e_i} e_i} J)(X)\\
			&=\sum_{i=1}^{2n} \nabla_{e_i} \nabla_{e_i} JX-2\nabla_{e_i} J \nabla_{e_i} X+J\nabla_{e_i} \nabla_{e_i} X-(\nabla_{\nabla_{e_i} e_i} J)(X)
		\end{align*}
		Therefore, 
		\begin{align*}
			[J,\nabla^*\nabla J](X)&=\sum_{i=1}^{2n} \{ J\nabla_{e_i} \nabla_{e_i} JX-2 J\nabla_{e_i} J \nabla_{e_i}  X-\nabla_{e_i} \nabla_{e_i} X-J(\nabla_{\nabla_{e_i} e_i} J)(X)\\
			&\quad \qquad -(-\nabla_{e_i} \nabla_{e_i} X-2\nabla_{e_i} J \nabla_{e_i} JX+J\nabla_{e_i} \nabla_{e_i} JX-(\nabla_{\nabla_{e_i} e_i} J)(JX))\}\\
			&=2 \sum_{i=1}^{2n} \left(\nabla_{e_i} J \nabla_{e_i} J X-J\nabla_{e_i} J\nabla_{e_i} X -J(\nabla_{\nabla_{e_i} e_i} J)X\right),
		\end{align*}
		where we have used $J(\nabla_X J)=-(\nabla_X J) J$.
	\end{proof}
	
	Using this lemma, we characterize in the following result the harmonic almost complex structures on a metric almost abelian Lie algebra.
	
	\begin{theorem}\label{theorem:main}
		Let $(J,\pint)$ be an almost Hermitian structure on an almost abelian Lie algebra $\g$. Then $J$ is harmonic, i.e. $[J,\nabla^*\nabla J]=0$, if and only if the following conditions are satisfied:
		\begin{enumerate}
			\item[\ri] $\mu \gamma+D_s \gamma-(\Tr S)\rho-J' D_a J' \rho=0$,
			\item[\rii] $D_a J' D_a J'-J'D_a J' D_a+(\Tr S) [D_a,J']J'=0$,
		\end{enumerate}
		with the notation from \eqref{eq: desc. L}.
	\end{theorem}
	
	\begin{proof}
		We are going to use the following expression of the Levi-Civita connection on $\g$ associated to $\pint$ determined in \cite{Mi}: 
		\begin{equation}\label{eq:nabla} \nabla_{e_0} e_0=0,\quad \nabla_{e_0} u=A u,\quad \nabla_{u} e_0=-S u,\quad \nabla_{u} v=\la Su, v\ra e_0,
		\end{equation}
		where $u,v\in\u$. Now, the third term in the right-hand side of \eqref{eq:conm-explicit} is straightforward to compute:
		\begin{equation}\label{eq:summand}
			\sum_{i=0}^{2n-1} J(\nabla_{\nabla_{e_i} e_i} J)(x)=(\Tr S) J(\nabla_{e_0} J)(x)=(\Tr S) J[A,J]x, \qquad x\in\g.
		\end{equation}
		
		For simplicity let us set $H:=\frac12 [J,\nabla^* \nabla J]$. For $x\in\a$, using \eqref{eq:nabla} and \eqref{eq:summand} we have
		\begin{align*}
			H(x)&=\sum_{i=0}^{2n-1} \left( \nabla_{e_i} J \nabla_{e_i} J X-J\nabla_{e_i} J\nabla_{e_i} X\right)-(\Tr S)J[A,J]x \\
			&=\nabla_{e_0} J AJx-J\nabla_{e_0} JAx+\sum_{i=1}^{2n-1} \left(\nabla_{e_i} (\la Se_i, Jx\ra e_1)-J\nabla_{e_i}(\la Se_i, x\ra e_1) \right)\\
			&\quad -(\Tr S)J[A,J]x\\
			&=(AJAJ-JAJA)x+\sum_{i=1}^{2n-1} \left(\la Se_i, e_1\ra \la Se_i, Jx\ra e_0-\la Se_i, e_1\ra \la Se_i, x\ra e_1\right)\\
			&\quad -(\Tr S) J[A,J]x \\
			&=(AJAJ-JAJA)x+\la Se_1, SJx\ra e_0-\la Se_1, Sx\ra e_1-(\Tr S)J[A,J]x.
		\end{align*}
		Next, we compute
		\begin{align*} 
			\la Hx, e_1\ra&=\la x, JAJA e_1\ra-\la S^2 e_1, x\ra+(\Tr S) \la x, [A,J] e_0\ra\\
			&=\la x, JAJ\rho-S(\mu e_1+\gamma)+(\Tr S) \rho\ra\\
			&=\la x, JAJ\rho-\mu \gamma-S\gamma-(\Tr S)\rho\ra\\
			&=-\la x, \mu\gamma+D_s \gamma-(\Tr S)\rho-J' D_a J'\rho\ra.
		\end{align*}     
		Since $H$ skew-commutes with $J$ and $\a$ is $J$-invariant, we get \[\la Hx, e_0\ra=-\la HJx,e_1\ra=\la J' x, \mu \gamma+D_s \gamma-(\Tr S) \rho-J' D_a J'\rho\ra.\]
		Now, for any $y\in\a$,
		\begin{align*}
			\la Hx, y \ra&=\la (AJAJ-JAJA-(\Tr S) J[A,J])x, y\ra\\
			&=\la (D_a J' D_a J'-J' D_a J' D_a+(\Tr S) [D_a,J']J')x,y\ra.
		\end{align*}
		Therefore, $Hx=0$ for all $x\in \a$ if and only if
		\[\mu \gamma+D_s \gamma-(\Tr S)\rho-J' D_a J'\rho=0,\quad
		D_a J' D_a J'-J'D_a J' D_a+(\Tr S) [D_a,J']J'=0,\] which are the conditions in the statement.
		
		Since $\nabla^*\nabla J$ is skew-symmetric  (because $\nabla_V J$ is skew-symmetric for any $V$), we have that also $H$ is skew-symmetric. This implies that $\la He_0, e_0\ra=0$. Moreover, $\la He_0, e_1\ra=\la HJ e_0, e_0\ra=-\la e_1, He_0\ra$, so $\la He_0,e_1\ra=0$. From this and $\la He_0, x\ra=-\la e_0,Hx\ra$ we obtain that $He_0=0$ if and only if condition $\ri$ holds. Similarly with $He_1$. This completes the proof.
	\end{proof}
	
	\medskip 
	
	As a consequence, we obtain:
	
	\begin{corollary}\label{corollary:unimod}
		With the same hypothesis as in Theorem \ref{theorem:main}, if $\g$ is unimodular then $J$ is harmonic if and only if
		\begin{enumerate}
			\item[\ri] $\mu \gamma+D_s \gamma-J'D_a J'\rho=0$,
			\item[\rii] $D_a J' D_a J'=J'D_a J' D_a$.
		\end{enumerate}
	\end{corollary}
	
	\begin{proof}
		We only have to note that $\Tr S=\Tr(\ad_{e_0})=0$.
	\end{proof}
	
	\medskip 
	
	In the next corollary, we consider the case when the orthogonal almost complex structure is integrable. It was shown in \cite{LRV} (see also \cite{AO}) that $J$ is integrable if and only if $w_0=0$ and $[D,J']=0$, which is equivalent to $[D_s,J']=[D_a,J']=0$. 
	
	\begin{corollary}\label{corollary:integrable}
		With the same hypothesis as in Theorem \ref{theorem:main}, if $J$ is integrable then $J$ is harmonic if and only if
		\[ Dv_0=(\Tr D) v_0. \]
	\end{corollary}
	
	\begin{proof}
		Suppose $J$ is integrable. Then, the condition $[D_a,J']=0$ implies automatically condition (ii) in Theorem \ref{theorem:main}.
		
		Since $w_0=0$, we have that $\gamma=\rho=\frac{1}{2}v_0$ and thus the first condition becomes \[ 
		\mu v_0+D_s v_0-(\Tr S)v_0+D_a v_0=0,\] which is equivalent to $Dv_0=(\Tr S-\mu)v_0=(\Tr D) v_0.$ 
	\end{proof}
	
	\medskip
	
	It follows from Corollary \ref{corollary:integrable} that if $\g$ is unimodular (i.e. $\mu+\Tr D=0$) and $J$ is integrable, then $J$ is harmonic if and only if \begin{equation}\label{eq:unimod int}
		Dv_0=-\mu v_0.
	\end{equation} 
	
	\smallskip
	
	It is easy to verify that, in this case, equation \eqref{eq:unimod int} is equivalent to $L^2 e_1=\mu^2 e_1$, since $Lv_0=Dv_0$.
	
	\medskip
	
	\section{Gray-Hervella classes on almost abelian Lie algebras}\label{section:GH}
	
	Let $(M,J,g)$ be a $2n$-dimensional almost Hermitian manifold with associated Levi-Civita connection $\nabla$, and Kähler form $\omega$ given by $\omega(X,Y)=g(JX,Y)$.
	
	We recall the Gray-Hervella classification of almost Hermitian structures \cite{GH}. Let $\mathcal W$ denote the space of $(0,3)$-tensor fields with the same symmetries as $\nabla \omega$. This space is a $\operatorname{U}(n)$-module which admits a decomposition 
	$\W=\W_1\oplus\W_2\oplus\W_3\oplus\W_4$,
	where $\W_i$ are the irreducible subspaces of the representation. Taking direct sums among them a total of 16 classes is obtained, and each one of these classes is defined by the vanishing of certain tensors. Some of these tensors are well known, for instance $d\omega$, $\delta\omega$ (where $\delta$ is the codifferential), the Nijenhuis tensor $N$ given in \eqref{eq:N}, and the Lee form $\theta$, which is the 1-form defined by 
	\begin{equation}\label{eq:Lee form}
		\theta(X)=-\frac{1}{n-1} \delta\omega(JX),\quad X\in\X(M).
	\end{equation} 
	
	Other tensors which characterize these classes are $T^{\pm}$ and $U$, given by 
	\begin{equation*} 
		T^{\pm}(X,Y,Z)=(\nabla_X \omega)(Y,Z)\pm(\nabla_{JX}\omega)(JY,Z),
	\end{equation*}
	\begin{equation}\label{eq:U}
		U(X,Y,Z)=g(X,Y) \delta\omega(Z)-g(X,Z)\delta\omega(Y)-g(X,JY) \delta\omega(JZ)+g(X,JZ)\delta\omega(JY).
	\end{equation}
	
	Then, the sixteen classes of almost Hermitian structures on manifolds of dimension $\geq6$ are shown in Table 1.
	
	\begin{table}[H]
		\centering
		\renewcommand*{\arraystretch}{1.5}
		\begin{tabular}{|c|c|}
			\hline
			Class & Defining conditions  \\ 
			\hline 
			$\{0\}$ & $\nabla \omega=0$ \\ 
			\hline 
			$\W_1$ & $3\nabla\omega=d\omega$ \\ 
			\hline  
			$\W_2$ & $d\omega=0$ \\ 
			\hline  
			$\W_3$ & $\delta\omega=0$, $N=0$  \\   
			\hline  
			$\W_4$ & $(\nabla_X\omega)(Y,Z)=-\frac{1}{2(n-1)} U(X,Y,Z) $\\
			\hline 
			$\W_1\oplus\W_2$ & $T^+=0$\\ 
			\hline 
			$\W_3\oplus \W_4$ & $N=0$ \\
			\hline 
			$\W_1\oplus\W_3$ & $T^-(X,X,Y)=0$, $\delta\omega=0$ \\
			\hline
			$\W_2\oplus\W_4$ & $d\omega=\theta\wedge \omega$ \\
			\hline 
			$\W_1\oplus\W_4$ & $(\nabla_X\omega)(X,Y)=-\frac{1}{2(n-1)} U(X,X,Y)$ \\
			\hline 
			$\W_2\oplus\W_3$ & $\sum\limits_{cyc} T^-(X,Y,Z)=0$, $\delta\omega=0$ \\
			\hline
			$\W_1\oplus\W_2\oplus\W_3$ & $\delta\omega=0$ \\
			\hline 
			$\W_1\oplus\W_2\oplus\W_4$ & $T^+(X,Y,Z)=-\frac{1}{n-1} U(X,Y,Z)$ \\
			\hline
			$\W_1\oplus\W_3\oplus\W_4$ & $\la N(X,Y), X\ra=0$ \\
			\hline
			$\W_2\oplus\W_3\oplus\W_4$ & $\sum\limits_{cyc}  T^-(X,Y,Z)=0$\\
			\hline 
			$\W$  & \text{No condition} \\ \hline
		\end{tabular}
		\caption{Almost Hermitian manifolds of dimension $\geq 6$}
	\end{table}
	
	In dimension 4, some of the sixteen classes collapse and only four of them remain, as shown in Table 1.
	
	\begin{table}[H]
		\renewcommand*{\arraystretch}{1.5}
		\centering
		\begin{tabular}{|c|c|}
			\hline
			Class & Defining conditions  \\ 
			\hline
			$\{0\}$ & $\nabla\omega=0$ \\
			\hline 
			$\W_2$ & $d\omega=0$ \\
			\hline
			$\W_4$ & $N=0$ \\
			\hline
			$\W$   & No condition \\ \hline
		\end{tabular}
		\caption{Almost Hermitian manifolds of dimension 4}
	\end{table}
	
	\medskip
	
	\begin{remark} \ 
		\begin{enumerate} 
			\item The manifolds belonging to some classes have special names, for instance: manifolds in class $\{0\}$ are Kähler; in $\W_1$,   nearly Kähler; in $\W_2$, almost Kähler; in $\W_3$, balanced; in $\W_1\oplus\W_2$, quasi-Kähler and in $\W_3\oplus \W_4$, Hermitian. In \cite{GH}, classes $\W_1\oplus\W_3\oplus\W_4$ and $\W_2\oplus\W_3\oplus\W_4$ are denoted $\mathcal{G}_1$ and $\mathcal{G}_2$, respectively.
			
			\item Locally conformally Kähler (LCK) manifolds, which are the Hermitian manifolds satisfying  $d\omega=\theta\wedge \omega$ and $d\theta=0$, belong to class $\W_4$. Moreover, when $\dim M\geq 6$ the class $\W_4$ coincides with the class of LCK manifolds. Indeed, if $J$ is in $\W_4$, using the general identity $d\omega(X,Y,Z)=(\nabla_X \omega)(Y,Z)+(\nabla_Y \omega)(Z,X)+(\nabla_Z \omega)(X,Y)$ for any $X,Y,Z$, we obtain
			\[ d\omega(X,Y,Z) =-\frac{1}{2(n-1)} [U(X,Y,Z)+U(Y,Z,X)+U(Z,X,Y)].\]
			Now, the cyclic sum $U(X,Y,Z)+U(Y,Z,X)+U(Z,X,Y)$ is equal to
			\begin{align*}
				&g(X,Y)\delta\omega(Z)-g(X,Z)\delta\omega(Y)-g(X,JY)\delta\omega(JZ)+g(X,JZ)\delta\omega(JY)\\
				& +g(Y,Z)\delta\omega(X)-g(Y,X)\delta\omega(Z)-g(Y,JZ)\delta\omega(JX)+g(Y,JX)\delta\omega(JZ)\\
				&+g(Z,X)\delta\omega(Y)-g(Z,Y)\delta\omega(X)-g(Z,JX)\delta\omega(JY)+g(Z,JY)\delta\omega(JX),
			\end{align*} and this is equal to
			\begin{align*}
				2[\omega(X,Y)\delta\omega(JZ)+\omega(Y,Z)\delta\omega(JX)+\omega(Z,X)\delta\omega(JY)].
			\end{align*}
			Therefore, using \eqref{eq:Lee form} we obtain that
			\[ d\omega=\theta \wedge \omega.\] Moreover, $d\omega=\theta\wedge\omega$ in dimension $\geq 6$ implies $d\theta=0$, so $J$ is LCK. Another consequence of this last implication (in $\dim\geq 6$) is that the fundamental 2-form $\omega$ of a manifold in class $\W_2\oplus\W_4$ is locally conformally symplectic \cite{V}. 
		\end{enumerate}
	\end{remark}
	
	\medskip 
	
	In the context of almost abelian Lie algebras, some Gray-Hervella classes have already been characterized. For instance, the class $\{0\}$ in \cite{FP-1}, $\W_2$ in \cite{LW}, $\W_3$ in \cite{FP}, LCK and $\W_2\oplus\W_4$ in \cite{AO}, $\W_3\oplus\W_4$ in \cite{LRV}, and it can be deduced from \cite{Bu} (see also \cite{AD}) that $\W_1$ collapses to $\{0\}$.
	
	Our aim in this section is to characterize all sixteen Gray-Hervella classes on almost abelian Lie algebras in a unified way. In the next section we will analyze the harmonicity condition in each of these classes. 
	
	To begin with, we need to describe the tensors appearing in Table 1 in terms of algebraic data.
	
	We start by computing the Nijenhuis tensor of $J$. Note that $N(x,x)=N(x,Jx)=0$ for all $x\in\g$. Also, $N(Jx,y)=-JN(x,y)$. Taking these symmetries into account and since the only non-zero brackets involve $e_0$, we can determine $N$ simply by computing $N(e_0,x)$, with $x\in\a$: 
	\begin{align}\label{eq: N en almost abelian}
		N(e_0,x)=(L+JLJ)x&=\la (L+JLJ)x,e_1\ra e_1+(D+J'DJ')x \nonumber\\
		&= \la (L^t+JL^tJ)e_1,x\ra e_1+(D+J'DJ')x \\
		&= \la w_0, x\ra e_1+(D+J'DJ')x. \nonumber 
	\end{align}
	
	In this invariant setting, the exterior derivative of $\omega$ can be computed using the formula \[d\omega(x,y,z)=-\omega([x,y],z)-\omega([y,z],x)-\omega([z,x],y).\]
	Since the only non-zero brackets involve $e_0$, we have, for $x,y\in \a$, $z\in\u$,
	\begin{align}\label{eq: d-omega} d\omega(e_0,e_1,x)&=-\la JLe_1, x\ra+\la JLx, e_1\ra=\la Le_1, Jx\ra=\la v_0, Jx\ra, \nonumber\\
		d\omega(e_0,x,y)&=-\la JLx, y\ra+\la JLy,x\ra=\la x, (D^t J'+J'D) y\ra, \\
		d\omega(z,x,y)&=0. \nonumber
	\end{align}
	
	Regarding the codifferential $\delta\omega$, its expression on an almost Hermitian Lie algebra was obtained for instance in \cite{AV}:
	\begin{equation*} 
		\delta\omega(x)=-\Tr \ad_{Jx}+\frac{1}{2} \la \sum_{i=0}^{2n-1} [Je_i, e_i], x\ra,
	\end{equation*}
	where $\{e_0,\ldots,e_{2n-1}\}$ is any orthonormal basis of $\g$. In the almost abelian context, since the only non-zero brackets involve $e_0$, this equation reduces to 
	\begin{equation*}
		\delta\omega(x)=-\Tr \ad_{Jx}-\la L e_1, x\ra.
	\end{equation*}
	Thus, 
	\begin{equation}\label{eq:delta omega} 
		\delta\omega(e_0)=0,\quad\delta\omega(e_1)=\Tr S-\mu=\Tr D,\qquad \delta\omega(x)=-\la v_0, x\ra, \; x\in \a.
	\end{equation}
	Note that 
	\begin{equation}\label{eq:delta=0}
		\delta\omega=0\iff \Tr D=0\quad \text{and}\quad v_0=0.
	\end{equation}
	
	It follows from \eqref{eq:delta omega} that the Lee form $\theta=-\frac{1}{n-1} (\delta\omega \circ J)$ is given by:
	\begin{equation}
		\theta(e_0)=-\frac{1}{n-1} \Tr D,\quad
		\theta(x)=\frac{1}{n-1} \la v_0, Jx\ra,\qquad x\in \u. 
	\end{equation}
	Finally, we compute $\nabla \omega$ using the formula 
	\[(\nabla_x \omega)(y,z)=-\omega(\nabla_x y, z)-\omega(y,\nabla_x z).\] 
	
	Note that, since $\nabla_x$ is skew-symmetric, we have 
	\begin{equation}\label{eq:J-symmetric}
		(\nabla_{x}\omega)(Jy,z)=(\nabla_x \omega)(y,Jz).
	\end{equation}
	In particular, $(\nabla_{x}\omega)(y,Jy)=0$ for all $x,y\in\g$. Hence, $\nabla\omega$ is given by
	\begin{align}\label{eq:nabla-omega}
		(\nabla_{e_0} \omega)(e_0,x)&=\la \rho, x\ra, \nonumber\\
		(\nabla_{e_0} \omega)(x,y)&=\la [A,J] x, y\ra, \nonumber\\
		(\nabla_{e_1} \omega)(e_1,x)&=\la \gamma, x\ra,  \\
		(\nabla_{z} \omega)(x,y)&=0, \nonumber\\
		(\nabla_{x} \omega)(e_1,y)&=\la Sx, y\ra \nonumber,
	\end{align} where $x,y\in\a$, $z\in \u$.
	The other cases follow from \eqref{eq:J-symmetric}.
	
	Next, we compute the aforementioned tensors $T^\pm$ and $U$. It follows from \eqref{eq:J-symmetric} that \begin{align*}   T^\pm(x,Jy,z)&=T^\pm(x,y,Jz),\\  T^\pm(x,y,z)&=-T^\pm(x,z,y). \end{align*} Moreover, it is easy to verify that $T^\pm(Jx,y,z)=\mp T^\pm(x,Jy,z)$. Then, $T^\pm$ is given by
	\begin{align}\label{eq:T-calc}
		T^\pm(e_0,e_0,x)&=\la \rho \pm \gamma, x\ra\nonumber,\\
		T^\pm(e_0,x,y)&=\la [A,J]x,y\ra,\\
		T^\pm(z,x,y)&=0\nonumber,\\
		T^\pm(x,e_1,y)&=\la (S\mp JSJ)x,y\ra \nonumber,  
	\end{align} where $x,y,z\in\a$. The other possibilities follow from the symmetries just described.
	
	On the other hand, using \eqref{eq:U} it is straightforward to verify that $U(x,z,y)=-U(x,y,z)$, $U(x,Jy,z)=U(x,y,Jz)$ and $U(Jx,Jy,z)=U(x,y,z)$. Thus, the tensor $U$ is determined by
	\begin{align*}
		U(e_0, e_0, x)&=-\la v_0, x\ra\nonumber,\\
		U(e_0, x,y)&=0\nonumber,\\
		U(x,e_1,y)&=-\la x,y\ra \Tr D,\\
		U(z,x,y)&=-\la z,x\ra \la v_0, y\ra+\la z,y\ra \la v_0, x\ra +\la z,Jx\ra \la v_0, Jy\ra-\la z,Jy\ra \la v_0, Jx\ra\nonumber,
	\end{align*}
	where $x,y,z\in\a$.
	
	We will apply all the expressions found above in the proof of the next theorem, which is the main result of this section. In the statement of the theorem we will use the following identifications: 
	\begin{align*} 
		\gl(n-1,\C)&\cong \{X\in \gl(\a)\mid [X,J']=0\}=\{X\in \gl(\a)\mid [X_a,J']=[X_s,J']=0\},\\  
		\u(n-1)&\cong \{X\in\gl(\a)\mid [X,J']=0, X^t=-X\}=\{X\in \gl(\a )\mid X_s=0, [X_a,J']=0\},\\   
		\sp(n-1,\R)&\cong \{X\in \gl(\a)\mid X^t J'+J'X=0\}=\{X\in \gl(\a)\mid [X_a,J']=0, X_s J'+J' X_s=0\}, \\
		\operatorname{Sym}(2n-2)&\cong \{X\in\gl(\a) \mid X=X_s\}. 
	\end{align*}
	where $X_a=\frac12 (X-X^t)$ and $X_s=\frac12 (X+X^t)$ are the skew-symmetric and symmetric parts of $X$, respectively. 
	\ 
	
	In the sequel we will use the notation \comillas{$J\in$} to denote that the almost Hermitian structure belongs to a certain class. For instance, $J\in \W_1$ means that the almost Hermitian manifold is nearly Kähler.
	
	\begin{theorem}\label{theorem:GH classes}
		Let $(J,\pint)$ be an almost Hermitian structure on a $2n$-dimensional almost abelian Lie algebra $\g=\R\ltimes_L\R^{2n-1}$, with $L$ given as in \eqref{eq:matrix L} and $n\geq 3$. Then we have the following relations between the Gray-Hervella classes: 
		\begin{itemize}
			\item $\W_1=\{0\}$,
			\item $\W_1\oplus\W_i=\W_i$, $2\leq i\leq 4$,
			\item $\W_1\oplus\W_2\oplus\W_4=\W_2\oplus \W_4$,
			\item $\W_1\oplus\W_3\oplus\W_4=\W_3\oplus\W_4 $.
		\end{itemize}
		Moreover, the defining conditions for each Gray-Hervella class are given in the following table:
		
		\begin{center}
			\renewcommand*{\arraystretch}{1.5}
			\begin{tabular}{|c|c|}
				\hline
				\textrm{Class} & \text{Conditions} \\
				\hline
				$\{0\}$ & $v_0=w_0=0$, $D\in \u(n-1)$  \\
				\hline
				$\W_2$ &  $v_0=0$, $D\in \sp(n-1,\R)$   \\
				\hline
				$\W_3$ &  $\Tr D=0$, $v_0=w_0=0$, $D\in \gl(n-1,\C)$  \\
				\hline
				$\W_4$ & $v_0=w_0=0$, $D-\left(\frac{\Tr D}{2(n-1)}\right)\I\in \u(n-1)$  \\
				\hline
				$\W_2\oplus\W_3$ & $\Tr D=0	, v_0=0, D\in \Sym(2n-2)\oplus\u(n-1)$  \\
				\hline
				$\W_1\oplus\W_2\oplus\W_3$ & $\Tr D=0$, $v_0=0$  \\
				\hline
				$\W_2\oplus\W_4$ & $v_0=0$, $D-\left(\frac{\Tr D}{2(n-1)}\right) \I\in \sp(n-1,\R)$\\
				\hline
				$\W_3\oplus\W_4$ & $w_0=0$, $D\in \gl(n-1,\C)$  \\
				\hline
				$\W_2\oplus\W_3\oplus\W_4$ & $D\in \Sym(2n-2)\oplus \u(n-1)$\\
				\hline
				$\W$ & No condition \\\hline
			\end{tabular}
		\end{center}
		
	\end{theorem}
	
	\medskip
	
	\begin{proof} 
		We start by analyzing the class $\W_2\oplus\W_3\oplus\W_4$ since this will help us characterize some other classes.
		
		$\bullet$ Class $\W_2\oplus \W_3\oplus \W_4$: The condition is $\sum\limits_{cyc} T^-(x,y,z)=0$. Note that this condition is invariant under permutations of $x,y,z$. Furthermore, if $y=x$ or $y=Jx$ then using the symmetries of $T^-$ the condition is automatically satisfied.
		
		Additionally, it follows from \eqref{eq:T-calc} that $T^-(z,x,y)=0$ for all $x,y,z\in\a$. Therefore, we only have to check the condition for $(e_0,x,y)$, with $x,y\in\a$:
		\begin{align*}
			0&=T^-(e_0,x,y)+T^-(x,y,e_0)+T^-(y,e_0,x)\\
			&=\la [A,J]x,y\ra+T^-(x,e_1,Jy)-T^-(Jy,e_1,x)\\
			&=\la [A,J]x,y\ra+\la (S+JSJ)x,Jy\ra-\la (S+JSJ)Jy,x\ra\\
			&=\la [A,J]x,y\ra \quad (\text{since}\, (S+JSJ)^t=S+JSJ)\\
			&=\la [D_a,J']x,y\ra.
		\end{align*}
		In conclusion, $J\in\W_2\oplus\W_3\oplus\W_4$ if and only if $[D_a,J']=0$. Therefore, $D_a\in \u(n-1)$.
		
		\medskip
		
		$\bullet$ Class $\{0\}$: As $\nabla\omega=0$, using  \eqref{eq:nabla-omega}, we have that $\rho=\gamma=0$, which implies $v_0=w_0=0$. Also, $[J',D_a]=0$ and $D_s=0$, which is equivalent to $D\in \u(n-1)$. This coincides with \cite[Lemma 3.6] {FP-1}.
		
		\medskip
		
		$\bullet$ Class $\W_1$: Let us compute $3\nabla\omega=d\omega$ in $(x,e_1,y)$, with $x,y\in \a$. Using \eqref{eq:nabla-omega} and \eqref{eq: d-omega} we  obtain $3\la Sx,y\ra=0$. Therefore $D_s=0$ and then, evaluating now $3\nabla\omega=d\omega$ in $(x,e_0,y)$ we have $3\la [D_a,J']x,y\ra=\la x,(-D_aJ'+J'D_a)y\ra$, which is equivalent to $[D_a,J']=0$. This implies $J\in\W_1\cap\W_1^\perp=\{0\}$, so it is Kähler. This follows also from \cite[Theorem 6.2]{AD}. 
		
		\medskip
		
		$\bullet$ Class $\W_2$: Using \eqref{eq: d-omega} and the condition $d\omega=0$, we obtain $v_0=0$ and $D^t J'+J'D=0$. This characterization was previously obtained in \cite[Proposition 4.1]{LW}.
		
		\medskip
		
		$\bullet$ Class $\W_3$: Using \eqref{eq:delta omega} and \eqref{eq:N} together with the conditions $\delta\omega=0$ and $N=0$, it follows easily that $\Tr D_s=0$, $v_0=w_0=0$ and $D+J'DJ'=0$. The latter condition is equivalent to $[D,J']=0$, which means $D\in\gl(n-1,\C)$. This result already appeared in \cite[Theorem 3.1]{FP}.
		
		\medskip
		
		$\bullet$ Class $\W_4$: The condition is $(\nabla_x \omega)(y,z)=-\frac{1}{2(n-1)} U(x,y,z)$ for all $x,y,z\in \g$.
		
		First, note that both tensors are skew-symmetric in the last two variables and, moreover evaluating them in $(x,Jy,z)$ is the same as evaluating them in $(x,y,Jz)$. 
		Therefore, we only have to check the condition for $(e_0,e_0,x)$, $(e_0,x,y)$, $(e_1,e_1,x)$, $(z,x,y)$, $(x,e_1,y)$, where $x,y\in\a$, $z\in\u$. 
		
		Using \eqref{eq:nabla-omega} and \eqref{eq:U} we have that, evaluating in $(e_0,e_0,x), x\in\a$,  
		\[\la \rho, x\ra=\frac{1}{2(n-1)}\la v_0, x\ra.\]
		For $(e_1,e_1,x), x\in\a$, we have
		\[ \la \gamma, x\ra=\frac{1}{2(n-1)} \la v_0, x\ra.\] 
		Both  equations are equivalent to $\rho=\gamma=\frac{1}{2(n-1)} v_0$, and since $n\geq 3$ this is equivalent to $v_0=w_0=0$. Since $v_0=0$, both tensors vanish evaluated in $(z,x,y)$ with $z\in\u, x,y\in\a$.
		
		Next, if we check  the condition for $(e_0,x,y)$ with $x,y\in\a$ we obtain $[D_a,J']=0$. Finally, for $(x,e_1,y)$ with $x,y\in \a$, the condition becomes 
		\[\la S x,y\ra=\frac{1}{2(n-1)}\la (\Tr D_s) x,y\ra,\] which is equivalent to $D_s=\frac{\Tr D_s}{2(n-1)} \I=\frac{\Tr D}{2(n-1)}\I$.
		
		It is straightforward to see that the conditions $[D_a, J']=0$ and  $D_s=\frac{\Tr D}{2(n-1)}\I$  hold if and only if $D-\left(\frac{\Tr D}{2(n-1)}\right)\I\in \u(n-1)$. This coincides with the conditions obtained in \cite[Theorem 3.3]{AO}.

		\medskip
		
		$\bullet$ Class $\mathcal{W}_1\oplus\mathcal{W}_2$: From $T^+\equiv0$ and \eqref{eq:T-calc}, we have that for $x,y\in\a$, $0=T^+(e_0,x,y)=\la [D_a,J']x,y\ra$. Therefore $[D_a,J']=0$ and thus $J\in \W_1^\perp$. Since also $J\in \W_1\oplus\W_2$ we have $\W_1\oplus\W_2=\W_2$.
		
		\medskip
		
		$\bullet$ Class $\mathcal{W}_3\oplus \mathcal{W}_4$: It follows from \eqref{eq: N en almost abelian} that the vanishing of $N$ is equivalent to $w_0=0$ and $[D,J']=0$, which coincides with \cite{LRV}, as mentioned before. 
		
		\medskip
		
		$\bullet$ Class $\W_1\oplus\W_3$: Note that the condition $T^-(x,x,z)=0$ for all $x,z\in\g$ is equivalent to $T^-(x,y,z)=-T^-(y,x,z)$ for all $x,y,z\in\g$. Choosing  $(x,e_1,y)$ with $x,y\in\a$ we have
		\[\la (D_s+J'D_sJ')x,y\ra=-\la (D_s+J'D_sJ')y,x\ra.\] 
		Since $D_s+J'D_sJ'$ is symmetric, we obtain $D_s+J'D_sJ'=0$. Now, using this, we obtain $T^-(e_1,x,y)=-T^-(y,x,e_1)=0$ for all $x,y\in \a$. Thus, \[\la [D_a,J']J'x,y\ra=T^-(e_0,J' x,y)=T^-(e_1,x,y)=0.\] This implies $[D_a,J']=0$, so $J\in \W_1^\perp$. In consequence, $\W_1\oplus\W_3=\W_3$. 
		
		\medskip
		
		$\bullet$ Class $\mathcal{W}_2\oplus\mathcal{W}_4$: The defining condition is $d\omega=\theta\wedge\omega$. We compute first $\theta\wedge\omega$:
		\[ \theta\wedge\omega(x,y,z)=\theta(x)\omega(y,z)+\theta(y)\omega(z,x)+\theta(z)\omega(x,y).\]
		For $x,y\in\a$, $z\in\u$, we have
		\begin{align}\label{eq: w2+w4}
			\theta\wedge\omega(e_0,e_1,x)&=\frac{1}{n-1}\la v_0, Jx\ra,\nonumber\\
			\theta\wedge\omega(e_0,x,y)&=\frac{1}{n-1} \la x, Jy\ra \Tr D_s,\\ 
			\theta\wedge\omega(z,x,y)&=\frac{1}{n-1}\left(\la Jz,x\ra \la v_0,Jy\ra+\la Jy,z\ra \la v_0,Jx\ra+\la Jx,y\ra\la v_0,Jz\ra\right). \nonumber
		\end{align}
		
		From $d\omega(e_0,e_1,x)=\theta\wedge\omega(e_0,e_1,x)$, $x\in\a$, using \eqref{eq: d-omega} and \eqref{eq: w2+w4}, we have
		$ v_0=\frac{1}{n-1} v_0$, and since $n\geq 3$ we get $v_0=0$. This automatically implies $\omega\wedge\theta(x,y,z)=0=d\omega(x,y,z)$, $x,y,z\in\a$.
		
		Now, evaluating $d\omega=\theta\wedge\omega$ in $(e_0,x,y)$ with $x,y\in\a$ we have $D^t J'+J'D=\frac{1}{n-1} (\Tr D_s) J'$, or equivalently 
		$D-\frac{\Tr D}{2(n-1)} \I\in \sp(n-1,\R)$. This coincides with \cite[Theorem 4.1]{AO}.
		
		\medskip
		
		$\bullet$ Class $\W_1\oplus\W_4$: The defining condition is $(\nabla_x\omega)(x,z)+\frac{1}{2(n-1)} U(x,x,z)=0$ for all $x,z\in\g$, which is equivalent to \[(\nabla_y\omega)(x,z)+\frac{1}{2(n-1)} U(y,x,z)=-( (\nabla_x\omega)(y,z)+\frac{1}{2(n-1)} U(x,y,z)).\] 
		
		Evaluating both sides of the expression in $(x,e_0,z)$ with $x,z\in\a$ we obtain
		\[ \la [A,J]x,z\ra=-\la Sx,Jz\ra+\frac{\Tr D_s}{2(n-1)} \la x,Jz\ra,\]
		which is equivalent to $J'[D_a,J']=-D_s+\frac{\Tr D_s}{2(n-1)} \I$. Since the left hand side is skew-symmetric and the right hand is symmetric, both must be zero. In particular, $[D_a,J']=0$, so $J\in\W_1^\perp$ and $\W_1\oplus\W_4=\W_4$.
		
		\medskip
		
		$\bullet$ Class $\mathcal{W}_1\oplus \mathcal{W}_2\oplus \mathcal{W}_3$: As we already noted in \eqref{eq:delta=0}, $\delta\omega=0$ is equivalent to  
		\[\Tr D=0,\quad v_0=0.\]
		
		$\bullet$ Class $\W_2\oplus\W_3$: The defining conditions are $\delta\omega=0$ and $\sum_{cyc} T^-(x,y,z)=0$, which are separately the defining conditions of $\W_1\oplus\W_2\oplus\W_3$ and $\W_2\oplus\W_3\oplus\W_4$. Therefore, $J\in \W_2\oplus\W_3$ if and only if
		\[ \Tr D=0, \quad v_0=0 \quad \text{and} \quad [D_a,J']=0.\]
		
		\medskip
		
		$\bullet$ Class $\W_1\oplus\W_2\oplus\W_4$: The defining condition is $T^+=-\frac{1}{n-1} U$. We evaluate in $(e_0,x,y)$, $x,y\in\a$, and, using \eqref{eq:T-calc} and \eqref{eq:U}, we obtain $[D_a,J']=0$. Therefore $J\in\W_1^\perp$ and thus
		$\W_1\oplus\W_2\oplus\W_4=\W_2\oplus\W_4$.
		
		\medskip 
		
		$\bullet$ Class $\W_1\oplus\W_3\oplus\W_4$: The defining condition is $\la N(x,y), x\ra=0$ for all $x,y\in\g$, which is equivalent to $\la N(x,y),z\ra=-\la N(z,y),x\ra$ for all $x,y,z\in\g$. Since $N(x,y)=0$ for $x,y\in\a$, it follows that $\la N(e_0,x),y\ra=-\la N(y,x),e_0\ra=0$.
		
		Now, using \eqref{eq: N en almost abelian} we obtain
		\begin{align*} 
			0&=\la N(e_0,x),y\ra,\\
			&=\la (D+J'DJ')x,y\ra, \\
			&=\la (D_a+J'D_aJ'+D_s+J'D_sJ')x,y\ra. 
		\end{align*}
		
		Since $D_a+J'D_a J'$ is skew-symmetric and $D_s+J'D_s J'$ is symmetric this condition means $[D_a,J']=0$ and $[D_s, J']=0$. In particular, $J\in\W_1^\perp$ so $\W_1\oplus\W_3\oplus\W_4=\W_3\oplus\W_4$.
	\end{proof}
	
	\medskip
		
	\
	
	The Gray-Hervella classes of almost Hermitian 4-dimensional almost abelian Lie algebras are easily obtained following the lines of the proof of Theorem \ref{theorem:GH classes} and using Table 2.
	
	\begin{corollary}\label{corollary:dim 4}
		Let $(J,\pint)$ be an almost Hermitian structure on a 4-dimensional almost abelian Lie algebra $\g=\R\ltimes_L \R^3$, with $L=\left[\begin{array}{c|cc} 
			\mu & r&s \\
			\hline 
			p & a&b\\
			q& c&d \end{array}\right]$. The defining conditions for each Gray-Hervella class are given in the following table: 
		\begin{center}
			\renewcommand*{\arraystretch}{1.5}
			\begin{tabular}{|c|c|}
				\hline
				Class & Conditions \\
				\hline
				$\{0\}$&  $p=q=r=s=0$, $a=d=0$, $b=-c$  \\
				\hline
				$\W_2$ & $p=q=0$, $a+d=0$ \\
				\hline
				$\W_4$& $r=s=0$, $a-d=0$, $b+c=0$  \\
				\hline
				$\W$ & No condition \\
				\hline
			\end{tabular}
		\end{center}
	\end{corollary}
	
	\medskip
	
	\begin{remark}
		It was proved in \cite{GM} that almost Hermitian structures in class $\W_3$ (balanced) or locally conformally Kähler structures are harmonic. This fact is also clear on almost abelian Lie algebras from our classification above.
	\end{remark}
	
	\medskip
	
	\section{Examples}\label{section:examples}
	In this section, we provide several explicit examples of harmonic and non harmonic almost complex structures on almost abelian Lie groups and some associated compact solvmanifolds. Throughout this section all Lie algebras are unimodular. 
	
	\subsection{Dimension 4}
	
	Let $\g=\R\ltimes_L \R^3$ be a 4-dimensional unimodular almost abelian Lie algebra with orthonormal basis $\{e_0,e_1,e_2,e_3\}$ such that $\R^3=\text{span}\{e_1,e_2,e_3\}$,  and orthogonal almost complex structure $J$ given by $Je_0=e_1$ and $Je_2=e_3$. In this basis the matrix $L$ is given by $L=\left[\begin{array}{c|cc} 
		\mu & r&s \\
		\hline 
		p & a&b\\
		q& c&d \end{array}\right]$. In order to analyze the harmonicity of the almost complex structure $J$ we express the conditions in Corollary \ref{corollary:unimod} in terms of $L$.
	
	Taking into account that any $2\times 2$ skew-symmetric matrix commutes with the $2\times 2$ orthogonal almost complex structure, it is clear that condition (ii) of Corollary \ref{corollary:unimod} automatically holds in dimension 4. Moreover,  condition (i) is equivalent to  
	\begin{equation}\label{eq:dim 4}
		bq+cs-d(p+r)=0, \quad cp+br-a(q+s)=0.
	\end{equation}
	Therefore, $J$ is harmonic if and only if \eqref{eq:dim 4} holds.
	
	\smallskip
	
	\begin{example}\label{example:non int} (Harmonic non integrable) For $L_0=\left[\begin{array}{c|cc}0&1&0\\ \hline 1&0&0\\0&0&0\end{array}\right]$, the almost complex structure $J$ on $\g=\R\ltimes_{L_0} \R^3$ is non integrable since $w_0^t=(1,0)\neq 0$. In fact, $J$ is in the general class $\W$. Furthermore, it follows from \eqref{eq:dim 4} that $J$ is harmonic.
		
		We show next that the corresponding simply connected Lie group $G$ admits a countable number of non-isomorphic lattices. Indeed, for any $m\in \N$, $m\geq 3$, let $t_m=\log(\frac{m+ \sqrt{m^2-4}}{2})$. Then, the matrix
		$\exp(t_m L_0)=\matriz{\cosh t_m&\sinh t_m&0\\ \sinh t_m&\cosh t_m&0\\0&0&1}$ is conjugate to the matrix $C_m=\begin{bmatrix} 0&-1&0\\1&m&0\\0&0&1\end{bmatrix}$, which has integer coefficients. It follows from Proposition \ref{proposition:bock} that any choice of $m$ gives rise to a lattice $\Gamma_m$ in $G$. According to Remark \ref{remark:gamma-isomorphic}, $\Gamma_m$ is isomorphic to $\tilde{\Gamma}_m=\Z\ltimes_{C_m} \Z^3$ and since $\tilde{\Gamma}_m/[\tilde{\Gamma}_m,\tilde{\Gamma}_m]\cong \Z^2\oplus \Z_{m-2}$, we have that $\Gamma_m$ is isomorphic to $\Gamma_n$ if and only if $m=n$. 
		
		The Lie algebra $\g$ is isomorphic to $\R \times \mathfrak{e}(1,1)$, where $\mathfrak{e}(1,1)$ is the Lie algebra of the group of rigid motions $E(1,1)$ of Minkowski 2-space, and $G$ is isomorphic to a product of $\R$ with the identity component $E_0(1,1)$ of $E(1,1)$, which is simply connected.
		
		\end{example}
		
		\medskip

		\begin{example} (Harmonic $\W_2$ - almost Kähler)
			According to \eqref{eq:dim 4} and Corollary \ref{corollary:dim 4} an almost Kähler structure is harmonic if and only if 
			\begin{equation}\label{eq:aK dim 4} 
				cs+ar=0, \quad br-as=0.
			\end{equation}
			Therefore, we have harmonic and non harmonic examples.
			
			For example, for $L_1=\left[\begin{array}{c|cc}
				0&&\\
				\hline 
				&1&\\&&-1 \end{array}
			\right]$ the almost complex structure $J$ is almost Kähler and harmonic, whereas for $L_2=\left[\begin{array}{c|cc} 0&1&0\\
				\hline 
				0&1&0\\0&0&-1\end{array}\right]$ the almost Kähler structure $J$ is not harmonic. We note that  the Lie algebras $\g_1=\R\ltimes_{L_1} \R^3$ and $\g_2=\R\ltimes_{L_2} \R^3$, and $\R\times\e(1,1)$ from Example \ref{example:non int} are all isomorphic, since $L_0, L_1$ and $L_2$ have the same Jordan form (Lemma \ref{lemma:isomorphic}).
			
			Therefore, any solvmanifold associated to $\R\times E_0(1,1)$ admits at least three almost Hermitian structures $(J_i,g_i)$ for $i=0,1,2$, such that $J_0$ and $J_1$ are harmonic with respect to $g_0$ and $g_1$, respectively, $J_0\in \W$ and  $J_1$ is almost Kähler, and $J_2$ is almost Kähler but not harmonic with respect to $g_2$. Moreover, it can be seen that $g_0$ and $g_1$ are isometric.
		\end{example}
		
		\medskip
		
		\begin{example} (Integrable harmonic)
			Given a 4-dimensional almost abelian metric Lie algebra with an orthogonal integrable almost complex structure $J$, we have that $J$ is harmonic if and only if $(J,\pint)$ is LCK. Indeed, the integrability condition ($a=d$, $c=-b$, $r=s=0$) plus the harmonicity given by \eqref{eq:dim 4} imply that $bq-ap=0$ and $bp+aq=0$. This means that, either $(a,b)=(0,0)$, or $(a,b)\neq (0,0)$ and $(p,q)=(0,0)$. In any case, $(J,\pint)$ is LCK by \cite{AO}. Conversely, recall that any LCK structure is harmonic by \cite{GM}. In \cite{AO}, lattices in 4-dimensional LCK almost abelian Lie algebras were constructed. The corresponding compact complex surfaces are primary Kodaira surfaces or Inoue surfaces of type $S^0$.
		\end{example}
		
		\medskip
		
		\begin{example} (Integrable non harmonic)
			Let $\g=\R\ltimes_L \R^3$, with $L=\left[\begin{array}{c|cc}
				0&0&0\\
				\hline 
				1&0&-1\\
				0&1&0
			\end{array}
			\right]$. Then, the almost complex structure $J$ above is integrable (\ref{corollary:dim 4}) but not harmonic since it does not satisfy \eqref{eq:dim 4}. 
			
			Note that the associated left invariant metric $g$ on the corresponding simply connected Lie group $G$ is not flat, since $L$ is not skew-symmetric (see \cite[Theorem 1.5]{Mi}). However, according to Remark \ref{remark:change}, if we consider the basis $\B'=\{e_0'=e_0+e_2, e_1'=e_1+e_3, e_2,e_3\}$ of $\g$, then the matrix $L$ takes the form $L'=\left[\begin{array}{c|cc} 0&0&0\\
				\hline 0&0&-1\\0&1&0\end{array}\right]$ and $J$ still satisfies $Je_0'=e_1'$ and $Je_2=e_3$. Let $\pint'$ be the inner product on $\g$ which makes $\B'$ orthonormal. Then, the corresponding left invariant metric $g'$ on $G$ is flat and moreover $(\g,J,g')$ is Kähler. In particular, $J$ is $g'$-harmonic. It is known that this Lie group $G$ admits lattices \cite{H}. Hence, for any lattice $\Gamma$ in $G$ the complex solvmanifold ($\Gamma\backslash G,J$) admits two invariant Hermitian metrics $g$ and $g'$ such that $J$ is not harmonic with respect to $g$, but it is harmonic with respect to $g'$ (moreover, $(J,g')$ is Kähler). 
			
			Next we analyze which solvmanifolds are obtained from $G$. It is well known that the only values of $t$ such that the matrix $\exp(tL)$ (or $\exp(tL')$) is conjugate to an integer matrix are $t\in \{2\pi, \pi, \frac{2\pi}{3}, \frac{\pi}{2}, \frac{\pi}{3}\}$ (see for instance \cite[Proposition 5.1]{Tol}). For $t=2\pi$, the corresponding lattice is abelian, isomorphic to $\Z^4$, thus the associated solvmanifold is diffeomorphic to the torus $\mathbb T^4$ (due to Corollary \ref{theorem:mostow}). For the other possible values of $t$, the associated solvmanifold is a hyperelliptic surface (see \cite{H}).  
		\end{example}
		
		\medskip 
		
		\begin{example}[Harmonic non integrable on a nilmanifold]
			
			\hfill
			
			(1)  Let $\g=\R\ltimes_L \R^3$ with $L=\left[\begin{array}{c|cc} 0&0&0\\ \hline 0&0&0\\0&1&0\end{array} \right]$. Then, it follows from Corollary \ref{corollary:dim 4} and \eqref{eq:aK dim 4} that $J$ is harmonic and almost Kähler. Moreover, $\g$ is the Lie algebra of the Lie group $H_3\times \R$, where $H_3$ is the 3-dimensional Heisenberg group. A nilmanifold associated to $H_3\times \R$ is the well-known Kodaira-Thurston manifold, and the invariant almost Kähler metric induced from $\g$ corresponds with the Abbena metric \cite{Abbena} (which was already seen to be harmonic in \cite{Wood}). 
			
			\smallskip
			
			(2) Let $\g=\R\ltimes_L\R^3$, where $L=\left[\begin{array}{c|cc} 0&1&0\\ \hline 0&0&0\\1&0&0\end{array}\right]$. According to \eqref{eq:dim 4}, the almost complex structure $J$ is harmonic and it is clearly non integrable; in fact, it belongs to the general class $\W$ (Corollary \ref{corollary:dim 4}). Since $L^3=0$ but $L^2\neq 0$, the Lie algebra $\g$ is 3-step nilpotent, and since the structure constants are rational, the corresponding Lie group $G$ admits lattices (\cite{Ma}). It is well known that $\g$ does not admit any complex structure (see for instance \cite{GR}). Moreover, using results of the theory of compact complex surfaces, a stronger result holds, namely any associated nilmanifold $\Gamma\backslash G$ does not admit complex structures, either invariant or not (\cite[Proposition 3.1]{BM}). 
			
			Furthermore, one can find almost Kähler structures on $\g$. Nevertheless, none of them are harmonic, due to the fact that if $L$ as in Corollary \ref{corollary:dim 4} satisfies $L^3=0$ and \eqref{eq:aK dim 4} then $L^2=0$.
		\end{example}
		
			\medskip 
		
		\subsection{Dimension \texorpdfstring{$\geq 6$}{}}
		
		We will provide examples of harmonic almost complex structures on $2n$-dimensional ($n\geq 3$) almost abelian unimodular Lie algebras $\g=\R\ltimes_L \R^{2n-1}$, belonging to classes $\W_2$, $\W_2\oplus\W_3$, $\W_1\oplus\W_2\oplus\W_3$, $\W_2\oplus\W_4$, $\W_3\oplus \W_4$, $\W_2\oplus\W_3\oplus \W_4$ and $\W$, respectively. We point out that the classes of the following examples of almost complex structures $J$ are genuine, in the sense that $J$ is not in any subclass contained in the class at issue. 
		
		In these examples we will always consider an orthonormal basis $\{e_0,e_1,\ldots,e_{2n-1}\}$ where $\u=\text{span}\{e_1,\ldots,e_{2n-1}\}$ and the orthogonal almost complex structure $J$ is given by $Je_{2i}=e_{2i+1}$ for $0\leq i\leq n-1$. All matrices will be expressed in this basis. Note that the unimodularity condition is 
		\begin{equation}\label{eq:mu-D}
			\mu+\Tr D=0,
		\end{equation} 
		with the notation from \eqref{eq:matrix L}. In each example we show the existence of lattices using Proposition \ref{proposition:bock}.
		
		\medskip
		
		\begin{example}
			
			(Harmonic $\W_2$ - almost Kähler)
			Let $\g=\R\ltimes_L \R^{2n-1}$ endowed with an almost Kähler structure $(J,\pint)$, that is, $J\in \mathcal{W}_2$. According to Theorem \ref{theorem:GH classes}, $D\in\sp(n-1,\R)$, therefore condition (ii) in Corollary \ref{corollary:unimod} is satisfied. Moreover $\Tr D=0$, and it follows from \eqref{eq:mu-D} that $\mu=0$. This together with condition (i) in Corollary \ref{corollary:unimod} imply that $J$ is harmonic if and only if \[D^t w_0=0.\]
			
			Let us choose\footnote{If $A$ and $B$ are matrices, $A\oplus B$ denotes the matrix $\matriz{A&0\\0&B}$. This naturally generalizes to the sum of $n$ matrices.} $D=\matriz{0&1\\1&0}^{\oplus (n-1)}\in \sp(n-1,\R)$, $\mu=0$ and $v_0=w_0=0$. It follows that $J$ is harmonic and $J\in \W_2$. Now, let $m\in \N, \, m\geq 3$, and 
			$t_m=\log(\frac{m+\sqrt{m^2-4}}{2})$. Then the matrix $\exp(t_m L)=(1)\oplus \matriz{\cosh t_m&\sinh t_m\\ \sinh t_m& \cosh t_m}^{\oplus (n-1)}$ is block-conjugate to $C_m=(1)\oplus \matriz{0&-1\\1&m}^{\oplus (n-1)}$. Applying Proposition \ref{proposition:bock} we obtain that for each $m$ there is a lattice $\Gamma_m$ and they are pairwise non-isomorphic since $\Gamma_m/[\Gamma_m,\Gamma_m]\cong \Z^2 \oplus (\Z_{m-2})^{n-1}$.
		\end{example}
			
		\medskip
		
		\begin{example} (Harmonic $\W_2\oplus\W_3$)
			Recall from Theorem \ref{theorem:GH classes} that $J\in\W_2\oplus\W_3$ if and only if $\mu=\Tr D=0$ (due to \eqref{eq:mu-D}), $v_0=0$ and $D_a J'=J' D_a$. Thus, condition (ii) from Corollary \eqref{corollary:unimod} is satisfied and condition (i) becomes simply
			\[ D^t w_0=0.\]   
			
			\smallskip
			
			Let $L=\left[\begin{array}{c|cccc}
				0&1&0&0&0\\ \hline
				0&0&0&0&0\\ 
				0&0&a_m&0&0\\
				0&0&0&-a_m&0\\
				0&0&0&0&0\end{array}
			\right]\oplus \matriz{a_m&0\\0&-a_m}^{\oplus (n-3)}$, where $a_m=\log(\frac{m+\sqrt{m^2-4}}{2})$, for $m\in \N$, $m\geq 3$. Clearly $J$ is harmonic and $J\in \W_2\oplus \W_3$. Furthermore, it is easy to verify that $J\notin \W_2\cup \W_3$. The corresponding simply connected Lie group admits lattices by Proposition \ref{proposition:bock} since \[\exp(L)=\matriz{1&1\\0&1}\oplus \operatorname{Diag}(e^{
				a_m},e^{-a_m},1)\oplus\matriz{e^{a_m}&0\\0&e^{-a_m}}^{\oplus (n-3)}\] is block-conjugate to an integer matrix. Therefore, we obtain solvmanifolds carrying a harmonic almost Hermitian structure which is genuinely $\W_2\oplus \W_3$.
		\end{example}
		
		\medskip
		
		\begin{example} (Harmonic $\W_1\oplus\W_2\oplus\W_3$) Let $J\in\W_1\oplus\W_2\oplus\W_3$, that is, $\mu=\Tr D$ and $v_0=0$. Since the Lie algebra is unimodular, it follows from \eqref{eq:mu-D} that $\mu=\Tr D=0$. Therefore, the harmonicity conditions from Corollary \ref{corollary:unimod} turn into  \[D_a J' D_a J'=J' D_a J' D_a, \quad D_s w_0+J'D_a J' w_0=0.\]
			
			For any $m\in \N, m\geq 3$, let $a_m=\log(\frac{m+\sqrt{m^2-4}}{2})$. A harmonic example is obtained by taking $\mu=0$, $v_0=w_0=0$ and  $D=\matriz{a_m&0&0&0\\0&0&-b&0\\0&b&0&0\\0&0&0&-a_m}\oplus \matriz{a_m&0\\0&-a_m}^{\oplus (n-3)}$,  for some $b\in\R$. Choosing $b\in\{2\pi, \pi, \frac{2\pi}{3}, \frac{\pi}{2}, \frac{\pi}{3}\}$ we have that $\exp(L)$ is block-conjugate to an integer matrix. Thus, using Proposition \ref{proposition:bock} we obtain  lattices and consequently solvmanifolds with a harmonic almost Hermitian structure $J$ in $\W_1\oplus\W_2\oplus\W_3$ and not in any subclass (since $[D_a,J']\neq 0$).
		\end{example}
		
		\medskip 
		
		\begin{example} (Harmonic $\W_2\oplus \W_4$) Let $J\in\W_2\oplus\W_4$, that is $v_0=0$ and $D=\lambda\I+B$, for some $\lambda\in\R$ where $B\in\sp(n-1,\R)$. Since $\g$ is unimodular, $\lambda=-\frac{\mu}{2n-2}$. We get that $[D_a,J']=[B_a,J']=0$, so
			condition (ii) in Corollary \ref{corollary:unimod} is fulfilled, and due to $v_0=0$ condition (i) becomes
			\[ D^t w_0=-\mu w_0.\]
			
			In \cite{AO}, examples of $2n$-dimensional almost abelian Lie algebras such that $J\in\W_2\oplus \W_4$ were given.  Explicitly, \[ L=(1)\oplus \operatorname{Diag}\left(0,\frac{1}{n},\ldots,\frac{n-1}{n},-\frac{1}{n},-\frac{2}{n},\ldots, -\frac{n-1}{n},-1\right).\] 
			$J$ is harmonic since $w_0=0$. It was also shown in \cite{AO} that there are values of $t\neq 0$ such that $\exp(tL)$ is conjugate to an integer matrix, and it was obtained a family $\{M_m\}_{m\in\N}$ of pairwise non homeomorphic solvmanifolds. 
			Therefore, these solvmanifolds carry a harmonic almost Hermitian structure in class $\W_2\oplus \W_4$. Note that $J\notin \W_2$ since $\mu\neq 0$ and $J\notin\W_4$ since $J$ is non integrable.
		\end{example}
		
		\medskip
		
		\begin{example} (Harmonic integrable) 
			
			Let $L_0=\left[\begin{array}{c|cccc}
				0&0&0&0&0\\
				\hline 
				0& 0&-a&0&0\\
				0& a&0 &0&0\\
				1&0&-a&0&0\\
				1&a&0&0&0
			\end{array}\right] \oplus \matriz{0&-b\\b&0}^{\oplus(n-3)}
			$, for some values of $a,b\in\R$. Then, the almost complex structure $J$ on the almost abelian Lie algebra $\g_0=\R\ltimes_{L_0} \R^{2n-1}$ is integrable due to Theorem \ref{theorem:GH classes}. Moreover, $J$ is harmonic since $Dv_0=0=\mu v_0$ (see \eqref{eq:unimod int}). Since $v_0\neq 0$, $J$ is neither balanced nor LCK, therefore $J\in\W_3\oplus\W_4$ but $J\notin \W_3\cup \W_4$. Furthermore, $L_0$ has the same Jordan form as \[L_1=\matriz{0&1\\0&0}\oplus\matriz{0&-a&0\\a&0&0\\0&0&0}\oplus\matriz{0&-b\\b&0}^{\oplus (n-3)}.\] Therefore, due to Lemma \ref{lemma:isomorphic}, $\g_0$ is isomorphic to the Lie algebra $\g_1=\R\ltimes_{L_1} \R^{2n-1}$. Choosing $a,b\in\{2\pi,\pi,\frac{\pi}{2},\frac{\pi}{3},\frac{2\pi}{3}\}$, $\exp(L_1)$ is block-conjugate to an integer matrix. Therefore the almost abelian simply connected Lie group $G$ admits lattices by Proposition \ref{proposition:bock} and thus, the corresponding solvmanifolds admit an integrable harmonic almost Hermitian structure.
		\end{example}
		
		\medskip
		
		\begin{example} (Harmonic $\W_2\oplus\W_3\oplus\W_4$)
			This class is characterized by $[D_a,J']=0$. Therefore  condition (ii) from Corollary \ref{corollary:unimod} is automatically satisfied and condition (i) becomes \[\mu\gamma+D_s \gamma-D_a \rho=0.\]
			For $m\in \N, m\geq 3$, let $a_m=\log(\frac{m+\sqrt{m^2-4}}{2})$.
			If $L=\matriz{0&0&2\\2&0&0\\0&0&0}\oplus \matriz{a_m&0\\0&-a_m}^{\oplus(n-2)}$, then $\mu=0$, $D_a=0$, $D_s=D$ and $\gamma=(1,1,0,\ldots,0)^t$. Consequently, $\mu\gamma+D_s\gamma-D_a\rho=0$ so $J$ is harmonic and clearly $J\in \W_2\oplus\W_3\oplus\W_4$. Since $v_0\neq 0$ and $w_0\neq 0$, $J$ is not in any subclass of $\W_2\oplus\W_3\oplus\W_4$. Furthermore, 
			\[\exp(L)=\matriz{1&0&2\\2&1&2\\0&0&1}\oplus\matriz{e^{a_m}&0\\0&e^{-a_m}}^{\oplus(n-2)}\] is block-conjugate to an integer matrix so by Proposition \ref{proposition:bock} we obtain lattices and hence solvmanifolds carrying a harmonic almost Hermitian structure $J\in\W_2\oplus\W_3\oplus\W_4$.  
		\end{example}
		
		\medskip
		
		\begin{example} (Harmonic $\W$)
			Let $L=\left[\begin{array}{c|cccc}
				0&0&0&0&2\\
				\hline 2&0&0&0&0\\0&0&0&-a&0\\0&0&a&0&0\\0&0&0&0&0
			\end{array}\right]\oplus\matriz{0&-b\\b&0}^{\oplus (n-3)}$ for some $a,b\in\R$.
			Then $\mu=0$, $\gamma=(1,0,0,1,0,\ldots,0)^t$, $\rho=(1,0,0,-1,0,\ldots,0)^t$, $D_s=0$ and $D_a=D$. Note that $[D_a,J']\neq 0$ but $D_a J' D_a J'=J' D_a J' D_a$ so condition (ii) in Corollary \ref{corollary:unimod} is satisfied. Since $\mu \gamma+D_s \gamma+J'D_aJ'\rho=0$, condition (i) is also satisfied and thus $J$ is harmonic. Given that $v_0\neq 0$, $w_0\neq 0$ and $[D_a,J']\neq 0$, we have that $J\in \W$ and not in any subclass. 
			
			The matrix $L$ has the same Jordan form as $L_1=\matriz{0&2&0\\0&0&2\\0&0&0}\oplus\matriz{0&-a\\a&0}\oplus\matriz{0&-b\\b&0}^{\oplus (n-3)}$, so  the Lie algebra $\g=\R\ltimes_L \R^{2n-1}$ is isomorphic to $\g_1=\R\ltimes_{L_1} \R^{2n-1}$ due to Lemma \ref{lemma:isomorphic}. For $a,b\in\{2\pi,\pi,\frac{2\pi}{3},\frac{\pi}{2},\frac{\pi}{3}\}$ we have that  \[\exp(L_1)=\matriz{1&2&2\\0&1&2\\0&0&1}\oplus\matriz{\cos a&-\sin a\\ \sin a & \cos a}\oplus\matriz{\cos b&-\sin b\\\sin b&\cos b}^{\oplus (n-3)}\] is block-conjugate to an integer matrix. According to Proposition \ref{proposition:bock} we obtain solvmanifolds carrying a generic harmonic almost Hermitian structure $J$.
		\end{example}
		
		\medskip
		
		\subsection{SKT structures and the harmonic condition}
		
		An interesting class of Hermitian manifolds, which does not fall into the Gray-Hervella classification, is given by \textit{strong K\"ahler with torsion} manifolds, or simply SKT manifolds (also known as pluriclosed). A Hermitian manifold $(M,J,g)$ is called SKT if $\partial\overline{\partial}\omega=0$, where $\omega=g(J\cdot,\cdot)$ denotes the fundamental $2$-form, as usual. These manifolds can be defined also by the equivalent condition $dc=0$, where $c(X,Y,Z):=g(T(X,Y),Z)$ is the torsion $3$-form of the Bismut connection on $M$ associated to the Hermitian structure $(J,g)$.
		
		Left-invariant SKT structures on almost abelian Lie groups have been studied in \cite{AL}, where the following result is proved\footnote{In \cite{AL} a different sign convention is used: there, $Je_0=-e_1$, but this does not affect the sign of $\mu$.}: 
		
		\begin{theorem}\label{theorem:SKT-teo}\cite[Theorem 4.6]{AL}
			Let $(J,\pint)$ be a Hermitian structure on a $2n$-dimen\-sional almost abelian Lie algebra $\g=\R e_0\ltimes_L \R^{2n-1}$, where $L$ is given by \eqref{eq:matrix L}. Then $(J,\pint)$ is SKT if and only if  $w_0=0$, $[D,J']=0$, $[D,D^t]=0$ and each eigenvalue of $D$ has real part equal to $0$ or $-\frac{\mu}{2}$. 
		\end{theorem}
		
		Using this characterization, we can determine when the complex structure of an SKT Hermitian structure is harmonic. As mentioned before, the condition $[D,J']=0$ is equivalent to $[D_s,J']=0$ and $[D_a,J']=0$, where $D_s$ (respectively, $D_a$) denotes the symmetric (respectively, skew-symmetric) part of $D$. On the other hand, the normality condition $[D,D^t]=0$ is equivalent to $[D_s,D_a]=0$. Therefore, $\{J',D_s,D_a\}$ is a commuting family of normal operators on $\mathfrak a$. Therefore, they are simultaneously unitarily diagonalizable over $\C$, so that there exists an orthonormal basis $\{ e_2,\ldots, e_{2n-1}\}$ of $\a$ in which $J'$ and $D$ have the following expressions:
		\begin{equation}\label{eq:J-D}
			J'=\matriz{0&-1\\1&0}^{\oplus(n-1)}, \quad 
			D=\matriz{a_1&-b_1\\b_1&a_1}\oplus\cdots \oplus \matriz{a_{n-1}&-b_{n-1}\\b_{n-1}&a_{n-1}},    
		\end{equation}   
		with $a_i=0$ or $a_i=-\dfrac{\mu}{2}$ and $b_i\in\R$, $1\leq i \leq n-1$, according to Theorem \ref{theorem:SKT-teo}. 
		
		\begin{proposition}\label{proposition:SKT-prop}
			With notation as above, if $(J,\pint)$ is SKT then $J$ is harmonic if and only if one of the following mutually exclusive conditions hold:
			\begin{enumerate}
				\item[$\ri$] $v_0=0$, or
				\item[$\rii$] $v_0\neq 0$ and $a_i=0$ for all $i$, that is, $D\in\u(n-1)$. Moreover, $D$ is singular and $v_0\in \Ker D$.
			\end{enumerate}
			Furthermore, if the almost abelian Lie algebra $\g$ is unimodular then:
			\begin{itemize}
				\item in Case $\ri$, either $\mu =0$ and $D\in\u(n-1)$, or $\mu\neq 0$ and there exists exactly one $i\in\{1,\ldots,n-1\}$ such that $a_i=-\frac{\mu}{2}$;
				\item in Case $\rii$, $\mu =0$.
			\end{itemize}
		\end{proposition}
		
		\begin{proof}
			Since $J$ is integrable, from Corollary \ref{corollary:integrable} we obtain that $J$ is harmonic if and only if $Dv_0=(\Tr D) v_0$.
			
			Clearly, if $v_0=0$ then $J$ is harmonic. Let us assume now that $J$ is harmonic and $v_0\neq 0$. This implies that $\Tr D$ is a real eigenvalue of $D$. According to Theorem \ref{theorem:SKT-teo}, we obtain $\Tr D=0$ or $\Tr D=-\dfrac{\mu}{2}$.
			
			Now, it follows from \eqref{eq:J-D} that $\Tr D=-k\mu$, where $k=|\{i\in\{1,\ldots,n-1\}\mid a_i=-\frac{\mu}{2}\}|$. Comparing with the possible values of $\Tr D$ (namely, $0$ or $-\frac{\mu}{2}$), we obtain that either $\mu=0$ or $\mu\neq 0$ and $k=0$. In any of these cases we obtain that $a_i=0$ for all $i$, so that $D\in \u(n-1)$. Moreover, $Dv_0=-k\mu v_0=0$.
			
			Finally, let us assume that $\g$ is unimodular, i.e. $\Tr L=0$. In Case $\ri$ we have $0=\Tr L=\mu+\Tr D=(1-k)\mu$, with $k$ as before. Therefore, either $\mu=0$ and $D\in \u(n-1)$, or $\mu \neq 0$ and $k=1$. In Case $\rii$ we have $0=\Tr L=\mu+\Tr D=\mu$ since $D\in \u(n-1)$.
		\end{proof}
		
		\smallskip 
		
		\begin{remark} \
			\begin{enumerate}
				\item We point out that in Case $\ri$ of Proposition \ref{proposition:SKT-prop}, with $\g$ unimodular and $\mu=0$, the Hermitian structure is actually K\"ahler, according to Theorem \ref{theorem:GH classes}. Moreover, the metric is flat, due to \cite[Theorem 1.5]{Mi}.
				
				\item SKT structures on almost abelian Lie algebras, and their relations with other geometric structures, were also studied in \cite{FP-1}. In particular, the authors obtain a classification, up to Lie algebra isomorphism, of the $6$-dimensional almost abelian Lie algebras that admit SKT structures.
			\end{enumerate}
		\end{remark}
		
		\smallskip
		
		\begin{corollary}
			For any $n\geq 2$, there exist $2n$-dimensional solvmanifolds admitting an SKT Hermitian structure $(J,g)$ such that $J$ is harmonic. 
		\end{corollary}
		
		\begin{proof}
			Consider the $2n$-dimensional almost abelian Lie algebra $\g=\R \ltimes_L \R^{2n-1}$ with the usual Hermitian structure $(J,\pint)$, where $L\in\gl(2n-1,\R)$ is given by
			\begin{equation}\label{eq:L-SKT}  L=
				\left[\begin{array}{c|cc}
					\mu &&\\ 
					\hline
					& -\frac{\mu}{2} &-a\\
					& a & -\frac{\mu}{2}
				\end{array}\right] \oplus \matriz{0&-b_1\\ b_1&0} \oplus \cdots\oplus \matriz{0&-b_{n-2} \\ b_{n-2} &0}, \quad \mu\neq 0.
			\end{equation}
			If $M(\mu,a)$ denotes the $(3\times 3)$-matrix appearing in \eqref{eq:L-SKT}, it was shown in \cite{AO} that for certain values of $(\mu,a), \, a\neq 0$, the matrix $\exp(M(\mu,a))$ is conjugate to an integer matrix. Therefore, choosing $b_i\in\{2\pi,\pi,\frac{2\pi}{3},\frac{\pi}{2},\frac{\pi}{3}\}$, we have that $\exp(L)$ is block-conjugate to an integer matrix. According to Proposition \ref{proposition:bock}, the corresponding simply connected almost abelian Lie group admits lattices and therefore the associated solvmanifold carries an invariant SKT structure with harmonic complex structure.
		\end{proof}
		
		\

		\ 
		

\begin{thebibliography}{99}
			
			\bibitem{Abbena}
			Abbena, E.: An example of an almost Kähler manifold which is not Kählerian. Boll. Unione Mat. Ital., VI. Ser. A \textbf{ 3}, 383--392 (1984)
			
			\bibitem{AD}
			Andrada, A., Dotti, I.G.: Conformal Killing-Yano 2-forms. Differential Geom. Appl. \textbf{58}, 103--119 (2018)
			
			\bibitem{AO}
			Andrada, A., Origlia, M.: Lattices in almost abelian Lie groups with locally conformal K\"ahler or symplectic structures. Manuscripta Math. \textbf{155}, 389--417 (2018)
			
			\bibitem{AV}
			Andrada, A., Villacampa, R.: Abelian balanced hermitian structures on unimodular Lie algebras. Transform. Groups \textbf{21}, 903--927 (2016)
			
			\bibitem{AT}
			Andrada, A., Tolcachier, A.: Harmonic complex structures and special Hermitian metrics on products of Sasakian manifolds. Preprint, arXiv:2301.09706 (2023)
			
			\bibitem{AL}
			Arroyo, R., Lafuente, R.: The long-time behavior of the homogeneous pluriclosed flow. Proc. Lond. Math. Soc. \textbf{119}, 266--289 (2019)
			
			\bibitem{BM}
			Bazzoni, G., Marrero, J.C.: Locally conformal symplectic nilmanifolds
			with no locally conformal K\"ahler metrics. Complex Manifolds \textbf{4}, 172--178 (2017)
			
			\bibitem{B} 
			Bock, C.: On low-dimensional solvmanifolds. Asian J. Math. \textbf{20}, 199--262 (2016)
			
			\bibitem{Bu}
			Butruille, J.-B.: Classification des
			variétés approximativement kähleriennes homogènes. Ann. Global Anal. Geom. \textbf{27}, 201–-225 (2005)
			
			\bibitem{DM1}
			Davidov, J., Mushkarov, O.: Harmonic almost-complex structures on twistor spaces. Israel J. Math. \textbf{131}, 319--332 (2002)
			
			\bibitem{DM2}
			Davidov, J., Mushkarov, O.: Harmonicity of the Atiyah-Hitchin-Singer and Eells-Salamon almost complex structures. Ann. Mat. Pura Appl. \textbf{197}, 185--209 (2018)
			
			\bibitem{DUM}
			Davidov, J., Ul Haq, A., Mushkarov, O.:
			Almost complex structures that are harmonic maps. J. Geom. Phys. \textbf{124}, 86--99 (2018)
			
			\bibitem{ES}
			Eells, J., Sampson, J.H.: Harmonic mapping Riemannian manifolds. Amer. J. Math. \textbf{86}, 109--160 (1964)
			
			\bibitem{FG}
			Fino, A., Grantcharov, G.: Properties of manifolds with skew-symmetric
			torsion and special holonomy. Adv. Math. \textbf{189}, 439--450 (2004)
			
			\bibitem{FKV}
			Fino, A., Kasuya, H., Vezzoni, L.: 
			SKT and tamed symplectic structures on solvmanifolds. Tohoku Math. J. \textbf{67}, 19--37 (2015)
					
			\bibitem{FP}
			Fino, A., Paradiso, F.: Balanced Hermitian structures on almost abelian Lie algebras. J. Pure Appl. Algebra \textbf{227}, Article 107186 (2022)
			
			\bibitem{FP-1}
			Fino, A., Paradiso, F.: Generalized K\"ahler almost abelian Lie groups. Ann. Mat. Pura Appl. \textbf{200}, 1781--1812 (2021)
			
			\bibitem{Fr}
			Freibert, M.: Cocalibrated structures on Lie algebras with a codimension one Abelian ideal. Ann. Glob. Anal. Geom. \textbf{42}, 537--563 (2012) 
			
			\bibitem{GM}
			Gonz\'alez-D\'avila, J.C., Mart\'in Cabrera, F.: Harmonic $G$-structures. Math. Proc. Camb. Phil. Soc. \textbf{146}, 435--459 (2009)
			
			\bibitem{GR}
			Goze, M., Remm, E.: Non existence of complex structures on filiform Lie algebras.
			Commun. Algebra \textbf{30}, 3777--3788 (2002) 
			
			\bibitem{GH}
			Gray, A., Hervella, L.: The sixteen classes of almost Hermitian manifolds and their linear invariants. Ann. Mat. Pura Appl. \textbf{123}, 35--58 (1980)
			
			\bibitem{H} 
			Hasegawa, K.: Complex and K\"ahler structures on compact solvmanifolds. J. Symplectic 
				Geom. \textbf{3}, 749--767 (2005)
			
			\bibitem{HL}
			He, W., Li, B.: The harmonic heat flow of almost complex structures. Trans. Amer. Math. Soc. \textbf{374}, 6179--6199 (2021)
			
			\bibitem{KL}
			Kath, I., Lauret, J.: A new example of a compact ERP $G_2$-structure. Bull. Lond. Math. Soc. \textbf{53}, 1692--1710 (2021) 
			
			\bibitem{LW}
			Lauret, J., Will, C.: On the symplectic curvature flow for locally homogeneous manifolds. J. Symplectic Geom. \textbf{15}, 1--49 (2017)
			
			\bibitem{LRV} 
			Lauret, J., Rodr\'iguez-Valencia, E.: On the Chern-Ricci flow and its solitons for Lie groups. Math. Nachr. \textbf{288}, 1512--1526 (2015)
			
			\bibitem{Ma} 
			Malcev, A.: On a class of homogeneous spaces. Izv. Akad. Nauk. Armyan. SSSR Ser.
				Mat. \textbf{13}, 9–32 (1949). English translation in: Amer. Math. Soc. Transl. 1951, No. \textbf{39}, 33 pp.(1951)
			
			\bibitem{Mi} 
			Milnor, J.: Curvatures of left invariant metrics on Lie groups. Adv. Math. \textbf{21}, 293--329 (1976)
			
			\bibitem{MSa}
			Moreno, A., Sá Earp, H.: Explicit soliton for the Laplacian co-flow on a solvmanifold. S\~ao Paulo J. Math. Sci. \textbf{15}, 280--292 (2021)
			
			\bibitem{Mo}
			Mostow, G.D.: Factor spaces of solvable groups. Ann. of Math. \textbf{60}, 1--27 (1954)
					
			\bibitem{Tol}
			Tolcachier, A.: Holonomy groups of compact flat solvmanifolds. Geom. Dedicata \textbf{209}, 95--117 (2020)
			
			\bibitem{V}
			Vaisman, I.: On locally conformal almost K\"ahler manifolds. Israel J. Math. \textbf{24}, 338--351 (1976)
			
			\bibitem{UV}
			Ugarte, L., Villacampa, R.: Symplectic harmonicity and generalized coeffective cohomologies. Ann. Mat. Pura Appl. \textbf{198}, 1351--1380 (2019)
			
			\bibitem{Wi} 
			Witte, D.: Superrigidity of lattices in solvable Lie groups. Invent. Math. {\bf 122}, 147--193 (1995)
			
			\bibitem{Wood-Crelle}
			Wood, C.M.: Instability of the nearly-K\"ahler six-sphere. J. Reine Angew. Math. \textbf{439}, 205--212 (1993)
			
			\bibitem{Wood}
			Wood, C.M.: Harmonic almost-complex structures. Compos. Math. \textbf{99}, 183--212 (1995)
		\end{thebibliography}
	\end{document}